\documentclass[12pt,a4paper,twoside]{article}

\usepackage{amsfonts,amsmath,amssymb,amsthm,latexsym}
\usepackage{times,mathptm}
\usepackage{color}
\usepackage{hyperref}
\usepackage[all]{xy}

\newtheorem{proposition}{Proposition}[section]
\newtheorem{lemma}[proposition]{Lemma}
\newtheorem{corollary}[proposition]{Corollary}
\newtheorem{theorem}[proposition]{Theorem}

\theoremstyle{definition}
\newtheorem{definition}[proposition]{Definition}
\newtheorem{example}[proposition]{Example}

\theoremstyle{remark}
\newtheorem{remark}[proposition]{Remark}

\newcommand{\thlabel}[1]{\label{th:#1}}
\newcommand{\thref}[1]{Theorem~\ref{th:#1}}
\newcommand{\selabel}[1]{\label{se:#1}}
\newcommand{\seref}[1]{Section~\ref{se:#1}}
\newcommand{\lelabel}[1]{\label{le:#1}}
\newcommand{\leref}[1]{Lemma~\ref{le:#1}}
\newcommand{\prlabel}[1]{\label{pr:#1}}
\newcommand{\prref}[1]{Proposition~\ref{pr:#1}}
\newcommand{\colabel}[1]{\label{co:#1}}
\newcommand{\coref}[1]{Corollary~\ref{co:#1}}
\newcommand{\relabel}[1]{\label{re:#1}}
\newcommand{\reref}[1]{Remark~\ref{re:#1}}
\newcommand{\exlabel}[1]{\label{ex:#1}}
\newcommand{\exref}[1]{Example~\ref{ex:#1}}
\newcommand{\delabel}[1]{\label{de:#1}}
\newcommand{\deref}[1]{Definition~\ref{de:#1}}
\newcommand{\eqlabel}[1]{\label{eq:#1}}
\newcommand{\equref}[1]{(\ref{eq:#1})}
\newcommand{\ot}{\otimes}

\newcommand{\cM}{\mathcal{M}}
\newcommand{\cC}{\mathcal C}
\newcommand{\cD}{\mathcal D}
\newcommand{\hH}{\mathcal H}
\newcommand{\can}{\mathrm{can}}
\def\stac#1{\raise-.2cm\hbox{$\stackrel{\,\displaystyle\otimes\,}{\scriptscriptstyle{#1}}$}}
\def\cstac#1{\raise-.3cm\hbox{$\stackrel{\,\displaystyle\Box\,}{\scriptscriptstyle{#1}}$}}
\def\smalldiamond{\raise.05cm\hbox{$\,{}_\lozenge$}}
\def\smallblackdiamond{\raise.05cm\hbox{$\,{}_\blacklozenge$}}

\renewcommand{\thefootnote}{\fnsymbol{footnote}}
\title{\MakeUppercase{Galois extensions over commutative and non-commutative
    base}} 
\author{Gabriella B\"ohm$^{1}$}
\date{}
\renewcommand{\date}{\vspace{-5mm}}
\begin{document}
\maketitle \vspace*{-3mm}\relax
\renewcommand{\thefootnote}{\arabic{footnote}}

\noindent \textit{\small $^1$ Research Institute for Particle and Nuclear
  Physics, Budapest, H-1525   Budapest 114, P.O.B.49,\\ 
Hungary\\
e-mail: G.Bohm@rmki.kfki.hu}

\begin{abstract}
\noindent
This paper is a written form of a talk. It gives a review of various notions
of Galois (and in particular cleft) extensions. Extensions by coalgebras,
bialgebras and Hopf algebras (over a commutative base ring) and by corings,
bialgebroids and Hopf algebroids (over a non-commutative base algebra) are
systematically recalled and compared.

\noindent
{\color{blue}
In the first version of this paper, the journal version of \cite[Theorem
2.6]{Brz:corext} was heavily used, in two respects. First, it was applied to 
establish an isomorphism between the comodule categories of two constituent
bialgebroids in a Hopf algebroid. Second, it was used to construct a Morita
context for any bicomodule for a coring extension. Regrettably, it turned out
that the proof of \cite[Theorem 2.6]{Brz:corext} contains an unjustified
step. Therefore, our derived results are not expected to hold at the stated
level of generality either.  
In the revised version we make the necessary corrections in both respects. 
In doing so, we obtain a corrected version of \cite[Theorem 4.2]{Bohm:hgdint}
as well, whose original proof contains a very similar error to  \cite[Theorem
2.6]{Brz:corext}.
}

\end{abstract}

\section*{Introduction}
The history of Hopf Galois extensions is nearly 40 years long, as it can be
traced back to \cite{ChaSwe:HopfGal}. Since then it is subject to a study of
always renewing interest. There are several reasons of this interest. First of
all, the algebraic structure is very rich. It has strong relations with the
problem of ring extensions. It is connected to (co-)module theory and a
descent problem. On the other hand, Hopf Galois theory unifies various
situations in an elegant manner. It is capable to describe e.g. classical
Galois extensions of fields or strongly group-graded algebras. Another
application of fundamental importance comes from non-commutative differential 
geometry. From this latter point of view, a (faithfully flat) Hopf Galois
extension is interpreted as a (dual form of a) non-commutative principal
bundle. 

Although the theory of Hopf Galois extensions was very fruitful, the
appearance of non-fitting examples forced it to be
generalized. Generalizations have been made in two different directions. In
one of them
the coacting Hopf algebra (or bialgebra) was replaced by a coalgebra. Since
this change results in loosing the monoidality of the category of comodules,
the notion of a comodule algebra is no longer available. Still, the essential
features of the theory turned out to be possible to maintain.
The notion of a Galois extension by a coalgebra appeared first in
\cite{Schn:PriHomSp} and then in \cite{BrzMaj:coabund}. The most general
definition, which will be used in this paper, can be found in the paper
\cite{BrzHaj:coalgGal} by Brzezi\'nski and Hajac.

In another direction of generalization, initiated by Kadison and
Szlach\'anyi in \cite{KadSzl:D2bgd}, the coacting bialgebra was replaced by a
bialgebroid over a non-commutative base algebra. Although in this case 
monoidality of the category of comodules persists, non-commutativity of
the base algebra results in a conceptually new situation. While modules over a
commutative ring form a {\em symmetrical} monoidal category, the monoidal
category of bimodules for any algebra is not symmetrical.

The two directions of generalization can be unified in the framework of
extensions by corings.

The aim of the current paper is to review all listed notions of Galois
extensions. We would like to show that -- after finding the proper,
categorically well established notions -- the theories over commutative and
non-commutative bases are reassuringly parallel. 

All Galois extensions in this paper are defined via {\em bijectivity} of a
certain canonical map. In the literature one can find a further
generalization, called a {\em weak} Galois extension, where the canonical map
is required to be only a {\em split monomorphism}
\cite[37.9]{BrzWis:coring}. Weak Galois extensions by coalgebras include
Galois extensions by weak Hopf algebras in \cite{CaeDeGr:WHA.Gal}. The issue
of weak Galois extensions is {\em not} considered in this paper.

The paper is organized as follows. In \seref{sec:prelims} we fix the notations,
and recall the basic notions, used later on. In \seref{sec:defs} 
definitions of the various Galois extensions are recalled. We start in
\seref{sec:Hopf.Gal} with the most classical and best understood case of a
Hopf Galois extension. Then in \seref{sec:coalg.Gal} we show how it fits the
more general case of a Galois extension by a coalgebra. These two sections
deal with definitions of Galois extensions over commutative base. In the case
of a non-commutative base, we proceed in a converse order. We start in
\seref{sec:coringGal} with recalling the most general instance of a Galois
extension by a coring. It is quite easily derived from the particular case of
a coalgebra. In \seref{sec:bgd.Gal} we define a Galois extension by
a bialgebroid as a special Galois extension by the underlying coring. 
Characterization is built on the monoidality of the category of comodules of a
bialgebroid. Finally, in \seref{sec:hgd.Gal} the coacting bialgebroid is
specialized to a Hopf algebroid. Here the complications are caused by the
presence of two bialgebroid structures, whose roles are clarified. All
sections are completed by examples.

An important (and relatively simple) class of Galois extensions is provided by
cleft extensions. In \seref{sec:cleft} the various notions of cleft extensions
are recalled. The order of the cases revisited follows the order in
\seref{sec:defs}. We present a unifying picture of {\color{blue}most} cleft
extensions, the one of cleft bicomodules, developed by Vercruysse and the
author in \cite{BohmVer:Mor&cleft}. It is shown that any cleft extension can
be characterized as a Galois extension with additional normal basis
property. Under mild further assumptions, a Strong Structure Theorem is proven
for cleft extensions. Cleft extensions by Hopf algebras or Hopf algebroids are
characterized as crossed products.

It has to be emphasized that this paper is of a review type. 
{\color{blue} Very few} new results
are presented. However, we hope that the way, they are collected from various
resources and re-organized here, give some new insight.
Since the results presented in the paper are already published, 
all proofs given here are sketchy. On the other hand, a big emphasis is put on
giving precise references to the original publications.

{\color{blue}
Parts of the first version of these notes (those concerning Galois extensions
with Hopf algebroids and those about cleft bicomodules) heavily relied on the
journal version of \cite[Theorem 2.6]{Brz:corext}. Regrettably, a few years
after publishing these notes it turned out that the proof of \cite[Theorem
2.6]{Brz:corext} contains an unjustified step. Although we are not aware of
any counterexamples, this gap forces us to revise our work, as we had to
revise references \cite{BohmBrz:hgd.cleft} and \cite{BohmVer:Mor&cleft}. Note 
similar errors also in \cite[Proposition 3.1]{Bohm:hgdGal} and \cite[Theorem
  4.2]{Bohm:hgdint}.  

There are two most important changes compared to the first version.
First, in the revised form of \seref{sec:hgd.Gal} we describe the Galois
theory of Hopf algebroids 
by distinguishing between comodules of the two constituent bialgebroids.
Second, we have to face the fact that the results in \cite{BohmVer:Mor&cleft}
are capable to handle cleft bicomodules only for so called {\em pure} coring
extensions. While this unifies cleft extensions by coalgebras (hence in
particular by Hopf algebras) and by corings, it is not known if all cleft
extensions by Hopf algebroids fit this framework. In
\seref{sec:nonpure.hgd.cleft} we propose a treatment of cleft extensions by
Hopf algebroids. Similarly to \cite{BohmVer:Mor&cleft}, a key role is played
by an appropriate Morita context. However, in general, this Morita context is
not known to be associated to a coring extension. 
Applying this Morita context,
for a cleft extension by an arbitrary Hopf algebroid, we prove a Strong
Structure Theorem (\thref{thm:nonpure.str.str.thm}). This yields in particular
a {\em corrected version} of \cite[Theorem 4.2]{Bohm:hgdint}, whose original
version is not known to be true because of an unjustified step in the proof,
see \reref{rem:fund.thm}. 

Using informality of the arXiv, as an experiment, for the convenience of the
readers of the first version, we write changes compared to the first version
in blue. 
}

\section{Preliminaries}\selabel{sec:prelims}

Throughout the paper we work over a commutative associative unital ring
$k$. 

The term {\bf $k$-algebra} (or sometimes only {\bf algebra}) means a
$k$-module $A$ equipped with a $k$-linear associative multiplication
$\mu:A\otimes_k A\to A$, with $k$-linear unit $\eta:k\to A$. On elements of
$A$ multiplication is denoted by juxtaposition. The unit element is denoted by
$1\in A$. We add labels, write $\mu_A$, $\eta_A$ or $1_A$, when the algebra
$A$ needs to be specified. 

The term {\bf $k$-coalgebra} (or sometimes only {\bf coalgebra}) means a
$k$-module $C$ equipped with a $k$-linear coassociative comultiplication
$\Delta:C\to 
C\otimes_k C$, with $k$-linear counit $\epsilon:C\to k$. On elements $c$ of
$C$, Sweedler's index notation is used for comultiplication:
$\Delta(c)=c_{(1)}\otimes_k c_{(2)}$ -- implicit summation understood.

A {\bf right module} of a $k$-algebra $A$ is a $k$-module $M$ equipped with a
$k$-linear associative and unital action $M\otimes_k A\to M$. On elements an
algebra action is denoted by juxtaposition. An {\bf $A$-module map} $M\to
M'$ is a $k$-module map which is compatible with the $A$-actions. The category
of right $A$-modules is denoted by $\cM_A$. For its hom sets the notation
$\mathrm{Hom}_A(-,-)$ is used. {\bf Left $A$-modules} are defined
symmetrically. Their category is denoted by ${}_A\cM$, and the hom sets by
${}_A\mathrm{Hom}(-,-)$. For two $k$-algebras $A$ and $B$, an $A$-$B$ {\bf
  bimodule} is a left $A$-module and right $B$-module $M$, such that the left
$A$-action $A\otimes_k M\to M$ is a right $B$-module map (equivalently, the
right $B$-action $M\otimes_k B\to M$ is a left $A$-module map). An $A$-$B$
{\bf bimodule map} is a map of left $A$-modules and right $B$-modules. The
category of $A$-$B$ bimodules is denoted by ${}_A\cM_B$ and the hom sets by
${}_A\mathrm{Hom}_B(-,-)$. 

A {\bf right comodule} of a $k$-coalgebra $C$ is a $k$-module $M$ equipped
with a $k$-linear coassociative and counital coaction $M\to M\otimes_k C$. On
elements a 
Sweedler type index notation is used for a coalgebra coaction: we write
$m\mapsto m_{[0]}\otimes_k m_{[1]}$ -- implicit summation understood. A {\bf
$C$-comodule map} $M\to M'$ is a $k$-module map which is compatible with the
$C$-coactions. The category of right $C$-comodules is denoted by $\cM^C$. For
its hom sets the notation $\mathrm{Hom}^C(-,-)$ is used. {\bf Left
$C$-comodules} are defined symmetrically. Their category is denoted by
${}^C\cM$, and the hom sets by ${}^C\mathrm{Hom}(-,-)$. For two $k$-coalgebras
$C$ and $D$, a $C$-$D$ {\bf bicomodule} is a left $C$-comodule and right
$D$-comodule $M$, such that the left $C$-coaction $M\to C\otimes_k M$ is a
right $D$-comodule map (equivalently, the right $D$-coaction $M\to M\otimes_k
D$ is a left $C$-comodule map). A $C$-$D$ {\bf bicomodule map} is a map of
left $C$-comodules and right $D$-comodules. The category of $C$-$D$ bicomodules
is denoted by ${}^C\cM^D$ and the hom sets by ${}^C\mathrm{Hom}^D(-,-)$.

We work also with bicomodules of a mixed type. For a $k$-algebra $A$ and a
$k$-coalgebra $C$, an $A$-$C$ {\bf bicomodule} is a left $A$-module and right
$C$-comodule $M$, such that the left $A$-action $A\otimes_k M\to M$ is a right
$C$-comodule map (equivalently, the right $C$-coaction $M\to M\otimes_k C$ is
a left $A$-module map). An $A$-$C$ {\bf bicomodule map} is a map of left 
$A$-modules and right $C$-comodules. The category of $A$-$C$ bicomodules is
denoted by ${}_A\cM^C$ and the hom sets by ${}_A\mathrm{Hom}^C(-,-)$.

For any algebra $R$, the category ${}_R\cM_R$ of $R$-$R$ bimodules is
monoidal. Monoidal product is given by the $R$-module tensor product, with
monoidal unit the regular bimodule. Hence the notion of a $k$-algebra (i.e. a
monoid in the monoidal category of $k$-modules) can be extended to an
arbitrary (non-commutative) base algebra $R$, as follows. An {\bf $R$-ring} is
a monoid in ${}_R \cM_R$. By definition, it means an $R$-$R$ bimodule $A$
equipped with an $R$-$R$ bilinear associative multiplication $\mu:A\otimes_R
A\to A$, with $R$-$R$ bilinear unit 
$\eta:R\to A$. Note that an $R$-ring $(A,\mu,\eta)$ can be characterized
equivalently by a $k$-algebra structure in $A$ and a $k$-algebra map
$\eta:R\to A$. A {\bf right module} of an $R$-ring $A$ is defined as an
algebra for the monad $-\otimes_R A:\cM_R\to \cM_R$. This notion coincides
with the one of a module for the respective $k$-algebra $A$.

Later in the paper we will be particularly interested in rings over a base
algebra $R\otimes_k R^{op}$, i.e. the tensor product of an algebra $R$
and its opposite $R^{op}$, which is the same $k$-module $R$ with opposite 
multiplication. In this case the unit map $\eta:R\otimes_k R^{op}\to A$ can be
equivalently given by its restrictions $s:=\eta(-\otimes_k 1_R):R\to A$ and
$t:=\eta(1_R\otimes_k -):R^{op}\to A$, with commuting ranges in $A$. The
algebra maps $s$ and $t$ are called the {\bf source} and {\bf target} maps,
respectively. Thus an $R\otimes_k R^{op}$-ring is given by a triple $(A,s,t)$,
consisting of a $k$-algebra $A$ and $k$-algebra maps $s:R\to A$ and
$t:R^{op}\to A$, with commuting ranges in $A$.

Just as the notion of a $k$-algebra can be generalized to an $R$-ring, also
that of a $k$-coalgebra can be generalized to an $R$-coring, for an arbitrary 
(non-commutative) base algebra $R$. By definition, an {\bf $R$-coring} is a
comonoid in ${}_R \cM_R$. It means an $R$-$R$ bimodule $\cC$ equipped with an
$R$-$R$ bilinear coassociative comultiplication $\Delta:\cC\to \cC\otimes_R
\cC$, with $R$-$R$ bilinear counit $\epsilon:\cC\to R$. Analogously to the
coalgebra case, on elements $c$ of $\cC$ Sweedler's index notation is used for
comultiplication: $\Delta(c)=c_{(1)}\otimes_R c_{(2)}$ -- implicit summation
understood. A {\bf right comodule} of an $R$-coring $\cC$ is defined as a
coalgebra for the comonad $-\otimes_R \cC:\cM_R\to \cM_R$. It means a right
$R$-module $M$ equipped with a right $R$-linear coassociative and counital
coaction $M\to M\otimes_R \cC$. On elements a Sweedler type index notation is
used for a coaction of a coring: we write $m\mapsto m_{[0]}\otimes_R m_{[1]}$
-- implicit summation understood. A {\bf right $\cC$-comodule map} $M\to M'$ is
a right $R$-module map which is compatible with the $\cC$-coactions. The
category of right $\cC$-comodules is denoted by $\cM^\cC$. For its hom sets
the notation $\mathrm{Hom}^\cC(-,-)$ is used. Left comodules and bicomodules
are defined analogously to the coalgebra case and also notations are analogous.

An element $g$ of an $R$-coring $\cC$ is said to be {\bf grouplike} if
$\Delta(g)=g\otimes_R g$ and $\epsilon(g)=1_R$. Recall a bijective
correspondence $g\mapsto (r\mapsto gr)$ between grouplike elements $g$ in
$\cC$ and right $\cC$-coactions in $R$, cf. \cite[Lemma 5.1]{Brz:str}.

A $k$-algebra map $\phi:R\to R'$ induces an $R$-$R$ bimodule structure in any
$R'$-$R'$ bimodule. What is more, it induces a canonical epimorphism
$\omega_\phi:M\otimes_R N\to M\otimes_{R'} N$, for any $R'$-$R'$ bimodules $M$
and $N$.
A {\bf map} from an $R$-ring $(A,\mu,\eta)$ to an $R'$-ring $(A',\mu',\eta')$
consists of a $k$-algebra map $\phi:R\to R'$ and an $R$-$R$ bimodule map
$\Phi:A\to A'$ (where the $R$-$R$ bimodule structure of $A'$ is induced by
$\phi$), such that 
$$
\Phi\circ \eta = \eta'\circ \phi,\qquad 
\Phi\circ \mu = \mu'\circ \omega_\phi\circ (\Phi\stac R \Phi).
$$
Note that $\Phi$ is necessarily a $k$-algebra map with respect to the
canonical $k$-algebra structures of $A$ and $A'$.

Dually, a {\bf map} from an $R$-coring $(\cC,\Delta,\epsilon)$ to an
$R'$-coring $(\cC',\Delta',\epsilon')$ consists of a $k$-algebra map $\phi:R\to
R'$ and an $R$-$R$ bimodule map $\Phi:\cC\to \cC'$ (where the $R$-$R$ bimodule
structure of $\cC'$ is induced by $\phi$), such that 
$$
\epsilon'\circ \Phi = \phi\circ \epsilon,\qquad 
\Delta'\circ \Phi = \omega_\phi\circ (\Phi\stac R \Phi)\circ \Delta.
$$
The notion of a coring extension was introduced in \cite[Definition
2.1]{Brz:corext}, as follows. An $R$-coring $\cD$ is a {\bf right extension}
of an $A$-coring $\cC$ if $\cC$ is a $\cC$-$\cD$ bicomodule, with left
$\cC$-coaction provided by the coproduct and some right $\cD$-coaction. 
{\color{blue}
Consider an $R$-coring $(\cD,\Delta_\cD,\varepsilon_\cD)$, which is a right
extension of an $A$-coring $(\cC,\Delta_\cC,\varepsilon_\cC)$. If the
equalizer  
\begin{equation}\eqlabel{eq:pure_eq}
\xymatrix{
M \ar[rr]^{\varrho}&&
M\stac A \cC \ar@<2pt>[rr]^-{\varrho\ot_A \cC}\ar@<-2pt>[rr]_-{M\ot_A
\Delta_\cC}&& 
M\stac A \cC\stac A \cC
}
\end{equation}
in $\cM_R$ is $\cD\ot_R \cD$-pure, i.e. it is preserved by the functor $-
\ot_R \cD \ot_R \cD :\cM_R \to \cM_R$, for any right $\cC$-comodule
$(M,\varrho)$, then we say that $\cD$ is a {\bf pure} coring extension of
$\cC$. By \cite[22.3]{BrzWis:coring} and its Erratum, by taking cotensor
products over $\cC$, 
a pure coring extension $\cD$ of $\cC$ induces a $k$-linear functor
$U=-\Box_\cC \cC:\cM^\cC\to \cM^\cD$ that commutes with the forgetful
functors $\cM^\cC\to \cM_k$ and $\cM^\cD\to \cM_k$, cf. the arXiv version of
\cite[Theorem 2.6]{Brz:corext}.
} 

\section{Definitions and examples}\selabel{sec:defs}

In this section various notions of Galois extensions are reviewed. We
start with the most classical notion of a Hopf Galois extension. The
definition is formulated in such a way which is appropriate for
generalizations. Generalizations are made in two directions. First, the
bialgebra symmetry in a Hopf Galois extension is weakened to
a coalgebra -- still over a commutative base. Then the base ring
is allowed to be non-commutative, so Galois extensions by corings are
introduced. It is understood then how the particular case of a Galois
extension by a bialgebroid is obtained. Finally we study Hopf algebroid Galois
extensions, first of all the roles of the two constituent bialgebroids.

\subsection{Hopf Galois extensions}\selabel{sec:Hopf.Gal}

Galois extensions of non-commutative algebras by a Hopf algebra have been
introduced in \cite{ChaSwe:HopfGal} and \cite{KreTak:HopfGal}, generalizing
Galois extensions of commutative rings by groups. Hopf Galois extensions unify
several structures (including strongly group graded algebras), studied
independently earlier. Beyond an algebraic importance, Hopf Galois extensions
are relevant also from the (non-commutative) geometric point of view. A Hopf
Galois extension (if it is faithfully flat) can be  interpreted as a (dual
version of a) non-commutative principal bundle.

\begin{definition}\delabel{def:bialg}
A {\em bialgebra} over a commutative ring $k$ is a $k$-module $H$, together
with a $k$-algebra structure $(H,\mu,\eta)$ and a $k$-coalgebra structure
$(H,\Delta,\epsilon)$ such that the counit $\epsilon: H\to k$ and the
coproduct $\Delta:H\to H\otimes_k H$ are $k$-algebra homomorphisms with
respect to the tensor product algebra structure of $H\otimes_k H$.
Equivalently, the unit $\eta:k\to H$ and the product $\mu:H\otimes_k
  H\to H$ are $k$-coalgebra homomorphisms with respect to the tensor product
  coalgebra structure of $H\otimes_k H$.

A {\em morphism of bialgebras} is an algebra and coalgebra map.

A bialgebra $H$ is a {\em Hopf algebra} if there exists a $k$-module map
$S:H\to H$, called the {\em antipode}, such that 
\begin{equation}\eqlabel{eq:Hopf.antip.ax}
\mu\circ (S\stac k H)\circ \Delta =\eta\circ \epsilon = \mu\circ(H\stac k
S)\circ \Delta.
\end{equation}
\end{definition} 

Since both the coproduct and the counit of a bialgebra $H$ are unital maps,
the ground ring $k$ is an $H$-comodule via the unit map $k\to H$. Furthermore,
since both the coproduct and the counit are multiplicative, the $k$-module
tensor product of two right $H$-comodules $M$ and $N$ is an $H$-comodule, with
the so called {\em diagonal coaction} $m\otimes_k n\mapsto m_{[0]}\otimes_k
n_{[0]}\otimes_k m_{[1]} n_{[1]}$. Since the coherence natural transformations
in $\cM_k$ turn out to be $H$-comodule maps with respect to these coactions,
the following theorem holds.

\begin{theorem}\thlabel{thm:bialg.com_cat_mon}
For a $k$-bialgebra $H$, the category of (left or right) comodules is a
monoidal category, with a strict monoidal forgetful functor to $\cM_k$.
\end{theorem}

\thref{thm:bialg.com_cat_mon} allows us to introduce a structure known
as an algebra extension by a bialgebra.

\begin{definition}\delabel{def:bialg.com_alg}
A right {\em comodule algebra} of a $k$-bialgebra $H$ is a monoid in the
monoidal category of right $H$-comodules. That is, an algebra and right
$H$-comodule $A$, whose multiplication and unit maps are right $H$-comodule
maps. Equivalently, $A$ is an algebra and right $H$-comodule such that the
coaction $\varrho^A:A\to A \otimes_k H$ is a $k$-algebra map, with respect to
the tensor product algebra structure in $A \otimes_k H$. 

The {\em coinvariants} of a right $H$-comodule algebra $A$ are the elements
of $A^{coH}\equiv \{\ b\in A\ |\ \varrho^A(ba)=b\varrho^A(a),\quad \forall
a\in A\ \}$. Clearly, $A^{coH}$ is a $k$-subalgebra of $A$. We say that $A$ is
a {\em (right) extension of} $A^{coH}$ by $H$.
\end{definition}

Note that, by the existence of a unit element $1_H$ in a $k$-bialgebra $H$,
coinvariants of a right comodule algebra $A$ (with coaction $\varrho^A$) can
be equivalently described as those elements $b\in A$ for which
$\varrho^A(b)=b\otimes_k 1_H$.

Let $H$ be a $k$-bialgebra, $A$ a right $H$-comodule algebra and
$B:=A^{coH}$. The right $H$-coaction $\varrho^A$ in $A$ can be used to
introduce a canonical map
\begin{equation}\eqlabel{eq:bialg.can}
\can:A\stac B A \to A \stac k H\qquad a\stac B a'\mapsto a\varrho^A(a').
\end{equation}

\begin{definition}\delabel{def:bialg.Gal}
A $k$-algebra extension $B\subseteq A$ by a bialgebra $H$ is said to be a {\em
  Hopf Galois extension} (or $H$-Galois extension) provided that the canonical
  map \equref{eq:bialg.can} is bijective.
\end{definition}

An important comment has to be made at this point. Although the structure
introduced in \deref{def:bialg.Gal} is called a {\em Hopf} Galois
extension $B\subseteq A$ (by a $k$-bialgebra $H$), the involved bialgebra $H$
is not required to be a Hopf algebra. However, in the most interesting case
when $A$ is a faithfully flat $k$-module, $H$ can be {\em proven} to be a
Hopf algebra, cf. \cite[Theorem]{Scha:Gal->Hopf}. 

\begin{example}\exlabel{ex:Hopf.Gal}
(1) Let $H$ be a bialgebra with coproduct $\Delta$ and counit $\epsilon$.
The algebra underlying $H$ is a right $H$-comodule algebra, with coaction
provided by $\Delta$. Coinvariants are multiples of the unit element. If $H$
is a Hopf algebra with antipode $S$, then the canonical map 
$$
\mathrm{can}:H\stac k H \to H\stac k H, \qquad h\stac k h'\mapsto
hh'_{(1)}\stac k h'_{(2)} 
$$
is bijective, with inverse $h\otimes_k h'\mapsto hS(h'_{(1)})\otimes_k
h'_{(2)}$. This proves that $k\subseteq H$ is an $H$-Galois
extension. Conversely, if the canonical map is bijective, then $H$ is a Hopf
algebra with antipode $(H\otimes_k \epsilon)\circ
\mathrm{can}^{-1}(1_H\otimes_k -)$. 

(2) {\em Galois extensions of fields.}
Let $G$ be a finite group acting by automorphisms on a field $F$. That is,
there is a group homomorphism $G\to \mathrm{Aut}(F)$, $g\mapsto \alpha_g$. For
any subfield $k$ of $F$, $F$ is a right comodule algebra of $k(G)$, the Hopf 
    algebra of $k$-linear functions on $G$. The coaction is given by $a\mapsto
    \sum_{g\in G} \alpha_g(a)\otimes_k \delta_g$, where the function $
    \delta_g\in k(G)$ takes the value $1_k$ on $g\in G$ and $0$ everywhere
    else. Denote $B:=F^{co k(G)}$. Then $F$ is a $k(G)$-Galois extension of
    $B$. That is, the canonical map 
$$
F\stac B F \to F\stac k k(G),\qquad a\stac B a'\mapsto \sum_{g\in G}
a\alpha_g(a')\stac k \delta_g
$$
is bijective \cite[Example 6.4.3]{DaNaRa:HopfAInt}.

(3) Let $A:=\oplus_{g\in G}\,\, A_g$ be a $k$-algebra graded by a finite group
    $G$. $A$ has a natural structure of a comodule algebra of the group Hopf
    algebra $kG$, with coaction induced by the map $a_g\mapsto a_g\otimes_k
    g$, on $a_g\in A_g$. The $kG$-coinvariants of $A$ are the elements of
    $A_{1_G}$, the component at the unit element $1_G$ of $G$. It is
    straightforward to see that the canonical map $A\otimes_{A_{1_G}} A \to A
    \otimes_k kG$ is bijective if and only if $A_g A_h=A_{gh}$, for all $g,h\in
    G$, that is $A$ is {\em strongly graded}.

\end{example}

\subsection{Galois extensions by coalgebras}\selabel{sec:coalg.Gal}

Comodules of any coalgebra $C$ over a commutative ring $k$
do not form a monoidal category. Hence one can not speak about comodule
algebras. Still, as it was observed in \cite{Schn:PriHomSp},
\cite{BrzMaj:coabund} and \cite{BrzHaj:coalgGal}, there is a
sensible notion of a Galois extension of algebras by a coalgebra.

\begin{definition}\delabel{def:coalg.ext}
Consider a $k$-coalgebra $C$ and a $k$-algebra $A$ which is a right
$C$-comodule via some coaction $\varrho^A:A\to A\otimes_k C$.
The {\em coinvariants} of $A$ are the elements of the $k$-subalgebra
$A^{coC}\equiv \{\ b\in A\ |\ \varrho^A(ba)=b\varrho^A(a),\quad \forall a\in A\
\}$. We say that $A$ is a {\em (right) extension of} $A^{coC}$ by $C$.  
\end{definition}

Let $B\subseteq A$ be an algebra extension by a coalgebra $C$.
One can use the same formula \equref{eq:bialg.can} to define a canonical map 
in terms of the $C$-coaction $\varrho^A$ in $A$,
\begin{equation}\eqlabel{eq:coalg.can}
\can:A\stac B A \to A \stac k C \qquad a\stac B a'\mapsto a\varrho^A(a').
\end{equation}

\begin{definition}\delabel{def:coalg.Gal}
A $k$-algebra extension $B\subseteq A$ by a coalgebra $C$ is said to be a {\em
  coalgebra Galois extension} (or $C$-Galois extension) provided that the
  canonical map \equref{eq:coalg.can} is bijective.
\end{definition}

\begin{remark}\relabel{rem:bialg<coalg.Gal}
In light of \deref{def:coalg.Gal}, a Hopf Galois extension $B\subseteq A$ by a
bialgebra $H$ in \deref{def:bialg.Gal} is the same as a Galois extension by
the coalgebra underlying $H$, such that in addition $A$ is a right
$H$-comodule algebra. 
\end{remark}

\begin{example}
Extending \cite[Example 3.6]{Schn:PriHomSp} (and thus \exref{ex:Hopf.Gal}
(1)), coalgebra Galois extensions can be constructed as in
\cite[34.2]{BrzWis:coring}. Examples of this class are called {\em quantum
homogenous spaces}.
Let $H$ be a Hopf algebra, 
with coproduct $\Delta$, counit $\epsilon$ and antipode $S$. Assume that $H$
is a flat
module of the base ring $k$. Let $A$ be subalgebra and a left coideal in
$H$. Denote 
by $A^+$ the augmentation ideal, i.e. the intersection of $A$ with the
kernel of $\epsilon$. The quotient of $H$ with respect to the coideal and
right ideal $A^+H$ is a coalgebra and right $H$-module ${\overline H}$. The
canonical epimorphism $\pi:H\to {\overline H}$ induces a right ${\overline
H}$-coaction $(H\otimes_k \pi)\circ \Delta$ on $H$. Denote the $\overline
H$-coinvariant subalgebra of $H$ by $B$. Clearly, $A\subseteq B$. Hence the
map 
$$
H\stac k {\overline H}\to H\stac B H,\qquad h\stac k \pi(h')\mapsto
hS(h'_{(1)})\stac B h'_{(2)}
$$
is well defined and yields the inverse of the canonical map
$$
H\stac B H\to H\stac k {\overline H}, \qquad h\stac B h'\mapsto hh'_{(1)}\stac k
\pi(h'_{(2)}). 
$$
This proves that $B\subseteq H$ is a Galois extension by ${\overline H}$.
Geometrically interesting examples of this kind are provided by Hopf
fibrations over the Podle\'s quantum spheres \cite{Brz:QHomSp},
\cite{MulSch:QHomSpff}. 
\end{example}

\subsection{Galois extensions by corings}\selabel{sec:coringGal}

If trying to understand what should be called a Galois extension by a coring,
one faces a similar situation as is the coalgebra case: comodules of an
arbitrary coring $\cC$ over an algebra $R$ do not form a monoidal category. As
the problem, also the answer is analogous.

\begin{definition}\delabel{def:coring.ext}
For a $k$- algebra $R$, consider an $R$-coring $\cC$ and an $R$-ring $A$ which
is a right $\cC$-comodule via some coaction $\varrho^A:A\to A\otimes_R 
\cC$. Require the right $R$-actions on $A$, corresponding to the
$R$-ring, and to the $\cC$-comodule structures, to be the same. 
The {\em coinvariants} of $A$ are the elements of the $k$-subalgebra
$A^{co\cC}\equiv \{\ b\in A\ |\ \varrho^A(ba)=b\varrho^A(a),\quad \forall a\in
A\ \}$. The $k$-algebra $A$ is said to be a {\em (right) extension of}
$A^{co\cC}$ by $\cC$.  
\end{definition}

It should be emphasized that though $A^{co\cC}$ in \deref{def:coring.ext} is a
$k$-subalgebra of $A$, it is not an $R$-subring (not even a one-sided
$R$-submodule) in general. 

Let $B\subseteq A$ be a $k$-algebra extension by an $R$-coring $\cC$.
Analogously to \equref{eq:coalg.can}, one defines a canonical map in terms of
the $\cC$-coaction $\varrho^A$ in $A$,
\begin{equation}\eqlabel{eq:coring.can}
\can:A\stac B A \to A \stac R \cC\qquad a\stac B a'\mapsto a\varrho^A(a').
\end{equation}

\begin{definition}\delabel{def:coring.Gal}
A $k$-algebra extension $B\subseteq A$ by an $R$-coring $\cC$ is said to be a
  {\em coring Galois extension} (or $\cC$-Galois extension) provided that the
  canonical map \equref{eq:coring.can} is bijective.
\end{definition}

\begin{example} 
In fact, in view of \cite[28.6]{BrzWis:coring}, any algebra extension
$B\subseteq A$, in which $A$ is a faithfully flat left or right $B$-module,
is a Galois extension by a coring, cf. a comment following
\cite[28.19]{BrzWis:coring}. Indeed, the obvious $A$-$A$ bimodule
$\cC:=A\otimes_B A$ is an $A$-coring, with coproduct and counit
$$
\Delta:\cC\to \cC\stac A \cC,\quad  a\stac B a'\mapsto (a\stac B 1_A)\stac A
(1_A\stac B a'),\qquad \textrm{and}\qquad
\epsilon:\cC\to A,\quad a\stac B a'\mapsto aa'.
$$
By the existence of a grouplike element $1_A\otimes_B 1_A$ in $\cC$, $A$ is a
right $\cC$-comodule via the coaction $A\to \cC$, $a\mapsto 1_A\otimes_B a$.
The coinvariants are those elements $b\in A$ for which $b\otimes_B
1_A=1_A\otimes_B b$. Hence, by the faithful flatness assumption made,
$A^{co\cC}=B$. The canonical map is then the identity map $A\otimes_B A$,
which is obviously bijective. 
\end{example}

\subsection{Bialgebroid Galois extensions}\selabel{sec:bgd.Gal}

In \reref{rem:bialg<coalg.Gal} we characterized Hopf Galois
extensions as Galois extensions $B\subseteq A$ by a constituent coalgebra in
a bialgebra $H$ such that $A$ is in addition an $H$-comodule algebra. The aim
of the current section is to obtain an analogous description of Galois
extensions by a bialgebroid $\hH$, replacing the bialgebra $H$ in 
\seref{sec:coalg.Gal}. In order to achieve this goal, a proper notion of a
`bialgebra over a non-commutative algebra $R$' is needed -- such that the
category of comodules is monoidal.
 
In \seref{sec:coringGal} we could easily repeat considerations in
\seref{sec:coalg.Gal} by replacing $k$-coalgebras -- i.e. comonoids in the
monoidal category $\cM_k$ of modules over a commutative ring $k$ -- by
$R$-corings --  i.e. comonoids in the monoidal category ${}_R\cM_R$ of
bimodules over a non-commutative algebra $R$. There is no such simple way to
generalize the notion of a bialgebra 
to a non-commutative base algebra $R$, by the following reason. While the
monoidal category $\cM_k$ is also {\em symmetrical}, not the bimodule category
${}_R\cM_R$. One can not consider bimonoids (bialgebras) in ${}_R\cM_R$. In
order to define a non-commutative base analogue of the notion of a bialgebra,
more sophisticated ideas are needed. The right definition was proposed in
\cite{Tak:Gr.alg} and independently in \cite{Lu:Hgd}. 

\begin{definition}\delabel{def:bgd}
A (right) {\em bialgebroid} over a $k$-algebra $R$ consists of an $R\otimes_k
R^{op}$-ring structure $(H,s,t)$ and an $R$-coring structure
$(H,\Delta,\epsilon)$ on the same $k$-module $H$, such that the following
compatibility axioms hold.
\begin{itemize}
\item[(i)] The $R$-$R$ bimodule structure of the $R$-coring $H$ is related to
  the $R\otimes_k R^{op}$-ring structure via
$$
rhr'=hs(r')t(r)\qquad \textrm{for } r,r'\in R\quad h\in H.
$$
\item[(ii)] The coproduct $\Delta$ corestricts to a morphism of $R\otimes_k
R^{op}$-rings 
$$
H\to H\times_R H\equiv \{\ \sum_i h_i\stac R h'_i\ |\ \sum_i 
s(r)h_i\stac R h'_i = \sum_i h_i\stac R t(r) h'_i,\quad \forall r\in R \ \},
$$
where $H\times_R H$ is an $R\otimes_k R^{op}$-ring via factorwise
multiplication and unit map $R\otimes_k R^{op}\to H\times_R H$, $r\otimes_k
r'\mapsto t(r') \otimes_R s(r)$.
\item[(iii)] The counit determines a morphism of $R\otimes_k R^{op}$-rings 
$$
H\to \mathrm{End}_k(R)^{op}, \qquad h\mapsto \epsilon\big(s(-)h\big),
$$
where $\mathrm{End}_k(R)^{op}$ is an $R\otimes_k R^{op}$-ring via
multiplication given by opposite composition of endomorphisms and unit map 
$R\otimes_k R^{op}\to\mathrm{End}_k(R)^{op}$, $r\otimes_k r'\mapsto r'(-)r$. 
\end{itemize}

A {\em morphism} from an $R$-bialgebroid $\hH$ to an $R'$-bialgebroid $\hH'$
is a pair consisting of a $k$-algebra map $\phi:R\to R'$ and an $R\otimes_k
R^{op}$-$R\otimes_k R^{op}$ bimodule map $\Phi:H\to H'$,
such that $(\phi\otimes_k \phi^{op}:R\otimes_k R^{op} \to R'\otimes_k R^{\prime
  op}, \Phi:H\to H')$ is a map from an $R\otimes_k R^{op}$-ring to an
$R'\otimes_k R^{\prime{op}}$-ring and $(\phi:R\to R',\Phi:H\to H')$ is a map
from an $R$-coring to an $R'$-coring. 
\end{definition}

A new feature of \deref{def:bgd}, compared to \deref{def:bialg}, is that here
  the monoid (ring) and comonoid (coring) structures are defined in {\em
  different} categories. This results in the quite involved form of the
  compatibility axioms. In particular, the $R$-module tensor product
  $H\otimes_R H$ is not a monoid in any category, hence the coproduct itself
  can not be required to be a ring homomorphism. It has to be corestricted to
  the so called {\em Takeuchi product} $H\times_R H$, which is an $R\otimes_k
  R^{op}$-ring, indeed. Similarly, it is not the counit itself which is a ring
  homomorphism, but the related map in axiom (iii).

Changing the comultiplication in a bialgebra to the opposite one (and leaving
the algebra structure unmodified) we obtain another bialgebra. If working over
a non-commutative base algebra $R$, one can replace  
the $R$-coring $(H,\Delta,\epsilon)$ in a right $R$-bialgebroid with the
co-opposite $R^{op}$-coring. Together with the $R^{op}\otimes_k R
$-ring $(H,t,s)$ (roles of $s$ and $t$ are interchanged!) they
form a right $R^{op}$-bialgebroid.

Similarly, changing the multiplication in a bialgebra to the opposite one (and
leaving the coalgebra structure unmodified) we obtain another
bialgebra. Obviously, \deref{def:bgd} is not invariant under the change of the 
$R\otimes_k R^{op}$-ring structure $(H,s,t)$ to $(H^{op},t,s)$ (and leaving the
coring structure unmodified). The $R\otimes_k R^{op}$-ring $(H^{op},t,s)$ and
the $R$-coring $(H,\Delta,\epsilon)$ satisfy symmetrical versions of the
axioms in \deref{def:bgd}. The structure in \deref{def:bgd} is usually termed
a {\em right $R$-bialgebroid} and the opposite structure is called a {\em left
  $R$-bialgebroid}. For more details we refer to \cite{KadSzl:D2bgd}. 

A most important feature of a bialgebroid for our application is formulated
in following \thref{thm:bgd.com_cat_mon}.
Comodules of an $R$-bialgebroid are meant to be comodules of the constituent
$R$-coring. \thref{thm:bgd.com_cat_mon} below was proven first in
\cite[Proposition 5.6]{Scha:bianc}, 
using an apparently more restrictive but in fact equivalent definition of a
comodule. The same version of \thref{thm:bgd.com_cat_mon} presented here can
be found in Section 2.2 of \cite{Bohm:hgdGal} or \cite[Proposition
1.1]{BalSzl:FinGal}.  

\begin{theorem}\thlabel{thm:bgd.com_cat_mon}
For a right $R$-bialgebroid $\hH$, the category of right comodules is a
monoidal category, with a strict monoidal forgetful functor to ${}_R\cM_R$. 
\end{theorem}

\begin{proof} (Sketch.)
Let $\hH$ be a right $R$-bialgebroid with $R\otimes_k R^{op}$-ring structure
$(H,s,t)$ and $R$-coring structure $(H,\Delta,\epsilon)$.
By the right $R$-linearity and unitality of the counit and the coproduct, the
base ring $R$ 
is a right comodule via the source map $s:R\to H$.
A right $\hH$-comodule $M$ (with coaction $m\mapsto m^{[0]}\otimes_R m^{[1]}$)
is a priori only a right $R$-module. Let us introduce a left $R$-action 
\begin{equation}\eqlabel{eq:left_R_mod}
rm := m^{[0]}\epsilon(s(r)m^{[1]}), \qquad \textrm{for } r\in R, m\in M.
\end{equation}
On checks that $M$ becomes an $R$-$R$ bimodule in this way, and any
$H$-comodule map becomes 
$R$-$R$ bilinear. Thus we have a forgetful functor $\cM^\hH\to {}_R\cM_R$. 
What is more, \equref{eq:left_R_mod} implies that, for any $m\in M$ and $r\in
R$, 
$$
rm^{[0]}\stac R m^{[1]}=m^{[0]}\stac R t(r) m^{[1]}.
$$
Hence the $R$-module tensor product of two right $\hH$-comodules $M$ and $N$
can be made a right $\hH$-comodule with the so called {\em diagonal coaction}
\begin{equation}\eqlabel{eq:prod_coac}
m\stac R n \mapsto m^{[0]}\stac R n^{[0]}\stac R m^{[1]}n^{[1]}.
\end{equation}
Coassociativity and counitality of the coaction \equref{eq:prod_coac} easily
follow by \deref{def:bgd}. 
The proof is completed by checking the right $\hH$-comodule map property of
the coherence natural transformations in ${}_R\cM_R$ with respect to the
coactions above.
\end{proof}

Let $\hH$ be a right bialgebroid. The category of right $\hH$-comodules, and
the category of left comodules for the co-opposite right bialgebroid
$\hH_{cop}$, are 
monoidally isomorphic. The category of right comodules for the opposite left
bialgebroid $\hH^{op}$ is anti-monoidally isomorphic to $\cM^\hH$. Thus (using
a bijective correspondence between left and right bialgebroid structures on a
$k$-module $H$, given by switching the order of multiplication), we conclude
by \thref{thm:bgd.com_cat_mon} that the category of right comodules of a left
$R$-bialgebroid is monoidal, with a strict monoidal forgetful functor to
${}_{R^{op}}\cM_{R^{op}}$.  

The important message of \thref{thm:bgd.com_cat_mon} is that there is a
sensible notion of a {\em comodule algebra} of a right bialgebroid
$\hH$. Namely, a 
right $\hH$-comodule algebra is a monoid in $\cM^\hH$. This means an $R$-ring 
and right $\hH$-comodule $A$ (with one and the same right $R$-module
structure), whose multiplication and unit maps are right 
$\hH$-comodule maps. Equivalently, $A$ is an $R$-ring and right $\hH$-comodule 
such that the coaction $\varrho^A$ corestricts to a map of $R$-rings
$$
A\to A\times_R H\equiv \{\ \sum_i a_i \stac R h_i \in A \stac R H\ |\ \sum_i
r a_i \stac R h_i =\sum_i a_i \stac R t(r)h_i,\quad \forall r\in R\ \}.
$$
The {\em Takeuchi product} $A\times_R H$ is an $R$-ring with
factorwise multiplication and unit map $r\mapsto 1_A \otimes_R s(r)$. Thus, in
analogy with \reref{rem:bialg<coalg.Gal}, we impose the following definition.

\begin{definition}\delabel{def:bgd.Gal}
A right {\em Galois extension $B\subseteq A$ by a right $R$-bialgebroid $\hH$}
is defined as a Galois extension by the $R$-coring underlying $\hH$, such that
in addition $A$ is a right $\hH$-comodule algebra. 
\end{definition}

\begin{remark}\relabel{rem:B_R}
Consider a right $R$-bialgebroid $\hH$, with structure maps denoted as in
\deref{def:bgd}, and a right $\hH$-comodule algebra $A$. Denote by $\eta:R\to
A$ the unit of the $R$-ring $A$.
By the left $B$-linearity and the unitality of the $\hH$-coaction on $A$, for
any element $b\in B:=A^{co\hH}$, $b^{[0]}\ot_R b^{[1]}=b\ot_R 1_H$.
By the right $R$-linearity and the unitality of the $\hH$-coaction on $A$, for
any element $r\in R$, $\eta(r)^{[0]}\ot_R \eta(r)^{[1]}=1_A\ot_R s(r)$.
Thus by the $\hH$-colinearity of the multiplication in $A$,
$$
(b\eta(r))^{[0]}\stac R (b\eta(r))^{[1]}=b\stac R s(r)=
(\eta(r) b)^{[0]}\stac R (\eta(r) b)^{[1]}.
$$
Applying $A\ot_R \varepsilon$ to both sides, we conclude that the subalgebra
$B$ of $A$ commutes with the range of $\eta$. 
\end{remark}

\begin{example}\exlabel{ex:depth2}
The {\em depth 2} property of an extension of arbitrary algebras $B\subseteq
A$ was introduced in \cite{KadSzl:D2bgd}, generalizing depth 2 extensions of
$C^*$-algebras. By definition, an algebra extension $B\subseteq A$ is right
depth 2 if there exists a finite integer $n$ such that $(A\otimes_B A)\oplus
-\cong \oplus^n A$, as {$A$-$B$} bimodules. By \cite[Theorem
5.2]{KadSzl:D2bgd} (see also \cite[Theorem 2.1]{Kad:Gal.bgd}), for a right
depth 2 algebra extension $B\subseteq A$, the centralizer $H$ of $B$ in the
obvious $B$-$B$ bimodule $A\otimes_B A$ has a right bialgebroid structure
$\hH$ over $R$, where $R$ is the commutant of $B$ in $A$. Its total algebra
$H$ is a finitely generated and projective left $R$-module. What is more, if
the algebra extension $B\subseteq A$ is also {\em balanced}, i.e. the
endomorphism algebra of $A$ as a left $\mathrm{End}_B(A)$-module is equal to
$B$, then $B\subseteq A$ is a right $\hH$-Galois extension. This example is a
very general one. By \cite[Theorem 2.1]{Kad:Gal.bgd} (see also 
\cite[Theorem 3.7]{BalSzl:FinGal}), any Galois extension by a finitely
generated and projective bialgebroid arises in this way.
\end{example}

After \deref{def:bialg.Gal} we recalled an observation in
\cite{Scha:Gal->Hopf} that a bialgebra, which admits a faithfully flat Hopf
Galois extension, is a Hopf algebra. Although the definition of a Hopf
algebroid will be presented only in forthcoming \seref{sec:hgd.Gal}, let us
anticipate here that no analogous result is known about bialgebroids. As a
matter of fact, the following was proven in \cite[Lemma
4.1.21]{Hobst:PhD}. Let $\hH$ be a right $R$-bialgebroid and $B\subseteq A$ a
right $\hH$-Galois extension such that $A$ is a faithfully flat left
$R$-module. Then the total algebra $H$ in $\hH$ is a right $\hH$-Galois
extension of the base algebra $R^{op}$ (via the target map). That is to say,
$\hH$ is a right {\em $\times_R$-Hopf algebra} in the terminology of
\cite{Sch:dua}. 
However, this fact does not seem to imply the existence of a Hopf algebroid
structure in $\hH$.

\subsection{Hopf algebroid Galois extensions}\selabel{sec:hgd.Gal}
 
As it is well known, the antipode of a $k$-Hopf algebra $H$ is a bialgebra map 
from $H$ to a Hopf algebra on the same $k$-module $H$, with opposite
multiplication and co-opposite comultiplication. We have seen in
\seref{sec:bgd.Gal} that 
the opposite multiplication in a bialgebroid does not satisfy the same axioms
the original product does. In fact, the opposite of a left bialgebroid is a
right bialgebroid and vice versa. Thus if the antipode in a Hopf algebroid
$\hH$ is expected to be a bialgebroid map between $\hH$ and a Hopf algebroid
with opposite multiplication and co-opposite comultiplication, then one has to
start with 
two bialgebroid structures in $\hH$, a left and a right one. The following
definition fulfilling this requirement was proposed in \cite{BohmSzl:hgdax},
where the antipode was required to be bijective. The definition was extended
by relaxing the requirement about bijectivity of the antipode in
\cite{Bohm:hgdint}. The set of axioms was reduced slightly further in
\cite[Remark 2.1]{BohmBrz:hgd.cleft}.

\begin{definition}\delabel{def:hgd}
A {\em Hopf algebroid} $\hH$ consists of a left bialgebroid $\hH_L$ over a
base algebra $L$ and a right bialgebroid $\hH_R$ over a base algebra $R$ on
the {\em same} total algebra $H$, together with a $k$-module map $S:H\to H$,
called the {\em antipode}. Denote the $L\otimes_k L^{op}$-ring structure in
$\hH_L$ by $(H,s_L,t_L)$ and the $L$-coring structure by $(H,\Delta_L,
\epsilon_L)$. Analogously, denote the $R\otimes_k R^{op}$-ring structure in
$\hH_R$ by $(H,s_R,t_R)$ and the $R$-coring structure by $(H,\Delta_R,
\epsilon_R)$. Denote the multiplication in $H$ (as an $L$-ring or as an
$R$-ring) by $\mu$. The compatibility axioms are the following.
\begin{itemize}
\item[(i)] The source and target maps satisfy the conditions
$$
s_L\circ \epsilon_L\circ t_R=t_R,\qquad
t_L\circ \epsilon_L\circ s_R=s_R,\qquad
s_R\circ \epsilon_R\circ t_L=t_L,\qquad
t_R\circ \epsilon_R\circ s_L=s_L.\\
$$
\item[(ii)] The two coproducts are compatible in the sense that 
$$
(\Delta_L\stac R H)\circ \Delta_R = (H\stac L \Delta_R)\circ \Delta_L,
\qquad  
(\Delta_R\stac L H)\circ \Delta_L = (H\stac R \Delta_L)\circ \Delta_R.
$$
\item[(iii)] The antipode is an $R$-$L$ bimodule map. That is, 
$$
S\big(t_L(l)ht_R(r)\big) = s_R(r) S(h) s_L(l),\qquad \textrm{for } r\in R,
l\in L, h\in H.
$$ 
\item[(iv)] The antipode axioms are
$$
\mu\circ(S\stac L H)\circ \Delta_L = s_R\circ \epsilon_R,\qquad
\mu\circ(H\stac R S)\circ \Delta_R = s_L\circ \epsilon_L.
$$
\end{itemize}
\end{definition}
Throughout these notes the structure maps of a Hopf algebroid $\hH$ will be
denoted as in \deref{def:hgd}. For the two coproducts $\Delta_L$ and
$\Delta_R$ we systematically use two versions of Sweedler's index notation:
we write 
$\Delta_L(h)=h_{(1)}\otimes_L h_{(2)}$ (with lower indices) and 
$\Delta_R(h)=h^{(1)}\otimes_R h^{(2)}$ (with upper indices), for $h\in H$.
In both cases implicit summation is understood.

The following consequences of the Hopf algebroid 
axioms in \deref{def:hgd} were
observed in \cite[Proposition 2.3]{Bohm:hgdint}. 

\begin{remark}
(1) The base algebras $L$ and $R$ are anti-isomorphic, via the 
map $\epsilon_R\circ s_L:L\to R$ (or $\epsilon_R\circ t_L:L\to R$), with
inverse $\epsilon_L\circ t_R:R\to L$ (or $\epsilon_L\circ s_R:R\to
L$). 

(2) The pair $(\epsilon_L\circ s_R:R\to L^{op}, S:H\to H^{op})$ is a morphism
of right bialgebroids, from $\hH_R$ to the opposite-co-opposite of $\hH_L$ and 
$(\epsilon_R\circ s_L:L\to R^{op}, S:H\to H^{op})$ is a morphism of left
bialgebroids, from $\hH_L$ to the opposite-co-opposite of $\hH_R$.
\end{remark}

{\color{blue}
Since a Hopf algebroid $\hH$ comprises two bialgebroid structures $\hH_L$ and
$\hH_R$, there are in general two different notions of their comodules. As it
turns out, the right definition of an $\hH$-comodule comprises both
structures. The following definition was proposed in \cite[Definition
  3.2]{Bohm:hgdGal} and \cite[Section 2.2]{BalSzl:FinGal}.

\begin{definition}\delabel{def:hgd_comod}
A {\em right comodule} of a Hopf algebroid $\hH$ is a right 
$L$-module as well as a right $R$-module $M$, together with a right
coaction $\varrho_R: M \to M\ot_R H$ of the constituent right bialgebroid
$\hH_R$ and a right coaction $\varrho_L:M\to M\ot_L H$ of the constituent
left bialgebroid $\hH_L$, 
such that $\varrho_R$ is an $\hH_L$-comodule map and $\varrho_L$ is an
$\hH_R$-comodule map. Explicitly, $\varrho_R$ is a right $L$-module map,
$\varrho_L$ is a right $R$-module map and
\begin{equation}\eqlabel{eq:hgd_comod}
(M\stac R \Delta_L)\circ \varrho_R = (\varrho_R \stac L H) \circ \varrho_L
\qquad \textrm{and}\qquad
(M\stac L \Delta_R)\circ \varrho_L = (\varrho_L \stac R H) \circ \varrho_R.
\end{equation}
Morphisms of $\hH$-comodules are meant to be $\hH_R$-comodule maps as well as
$\hH_L$-comodule maps. The category of right $\hH$-comodules is denoted by
$\cM^\hH$. 
\end{definition}
}
In the sequel we fix the following notation. For a Hopf
algebroid $\hH$, with constituent right bialgebroid $\hH_R$ and left
bialgebroid $\hH_L$, and a right $\hH$-comodule $M$, for $m\in M$ we write 
$m\mapsto m^{[0]}\otimes_R m^{[1]}$ and $m\mapsto m_{[0]}\otimes_L m_{[1]}$
for the $\hH_R$-, and $\hH_L$-coactions related by \equref{eq:hgd_comod}.
In both cases implicit summation is understood.

{\color{blue}
\begin{proposition}\prlabel{prop:ff}
For a Hopf algebroid $\hH$, with constituent right bialgebroid $\hH_R$ and left
bialgebroid $\hH_L$, the forgetful functor $\cM^\hH\to \cM^{\hH_R}$ is fully
faithful.
\end{proposition}
\begin{proof}
The forgetful functor $\cM^\hH\to \cM^{\hH_R}$ is obviously faithful. In order
to see that it is also full, consider two $\hH$-comodules $M$ and $M'$, and an
$\hH_R$-comodule map $f:M\to M'$. That is, a right $R$-module map $f$, such
that, for all $m\in M$,
\begin{equation}\eqlabel{eq:f_HR}
f(m)^{[0]}\stac R f(m)^{[1]}=f(m^{[0]})\stac R m^{[1]}.
\end{equation}
By definition, $M$ is a right $L$-module and by the right $L$-linearity of the 
$\hH_R$-coaction, 
\break
$ml=m^{[0]}\varepsilon_R (t_L(l)m^{[1]})$. Similarly, $M'$
is a right $L$-module and $f$ is clearly right $L$-linear. Regard $H$ as a left
$L$-module via 
$s_L$ and as a left $R$-module via $t_R$. Consider the well defined map
\begin{equation}\eqlabel{eq:Phi}
\Phi_{M'}:M'\stac R H \to M'\stac L H,\qquad 
m'\stac R h \mapsto m'_{[0]}\stac L m'_{[1]} S(h).
\end{equation}
Applying $\Phi_{M'}$ to both sides of \equref{eq:f_HR}
and using the right $\hH_R$-colinearity of the $\hH_L$-coaction on $M'$ and
one of the antipode axioms to simplify the left hand side, we obtain the
identity  
\begin{equation}\eqlabel{eq:f_2}
f(m)\stac L 1_H = f(m^{[0]})_{[0]} \stac L f(m^{[0]})_{[1]}S(m^{[1]}),
\end{equation}
for all $m\in M$. By the right $L$-linearity of $f$, \equref{eq:f_2} can be
used to compute
\begin{eqnarray*}
f(m_{[0]})\stac L m_{[1]} 
&=& f({m_{[0]}}^{[0]})_{[0]} \stac L f({m_{[0]}}^{[0]})_{[1]}S({m_{[0]}}^{[1]}) 
m_{[1]}\\
&=&  f(m^{[0]})_{[0]} \stac L f(m^{[0]})_{[1]}S({m^{[1]}}_{(1)}) {m^{[1]}}_{(2)}\\
&=&  f(m^{[0]})_{[0]} \stac L f(m^{[0]})_{[1]}
s_R\big(\varepsilon_R(m^{[1]}) \big)\\
&=&  \big(f(m^{[0]})\varepsilon_R(m^{[1]})\big)_{[0]} \stac L 
\big(f(m^{[0]})\varepsilon_R(m^{[1]})\big)_{[1]} 
=f(m)_{[0]}\stac L f(m)_{[1]},
\end{eqnarray*}
for all $m\in M$. 
In the first equality we applied \equref{eq:f_2}.
The second equality follows by the right $\hH_L$-colinearity of the
$\hH_R$-coaction on $M$.
In the third equality we used one of the antipode axioms in a Hopf algebroid. 
The penultimate equality follows by the right $R$-linearity of the
$\hH_L$-coaction on $M'$. 
The last equality follows by the right $R$-linearity of $f$ and the
counitality of the $\hH_R$-coaction on $M$.
This proves that $f$ is right $\hH_L$-colinear, hence it is
a morphism in $\mathrm{Hom}^\hH(M,M')$, as stated.
\end{proof}

As a simple consequence of \prref{prop:ff}, we obtain \cite[{\em
Corrigendum}, Proposition 3]{BohmBrz:hgd.cleft}.

\begin{corollary}\colabel{cor:hgd.coinv}
Let $\hH$ be a Hopf algebroid and $M$ be
a right $\hH$-comodule. Then any coinvariant of the $\hH_R$-comodule
$M$ is coinvariant also for the $\hH_L$-comodule $M$.
If moreover the antipode of $\hH$ is bijective then coinvariants of the
$\hH_R$-comodule $M$ and the $\hH_L$-comodule $M$ coincide. 
\end{corollary}

\begin{proof}
The base algebra $R$ of the constituent right bialgebroid $\hH_R$ is a right
$\hH_R$-comodule via the right regular $R$-action and the coaction $R\to
R\ot_R H\cong H$, $r\mapsto 1_R \ot_R s_R(r)\cong s_R(r)$ (determined by the
grouplike element $1_H$). 
Moreover, $R$ is also a right comodule of the constituent left $L$-bialgebroid
$\hH_L$, with right $L$-action $R\ot L\to R$, $r\ot l\mapsto
\varepsilon_R\big(t_L(l)\big) r$ and coaction $R \to R\ot_L H\cong H$,
$r\mapsto 1_R \ot_L s_R(r)\cong s_R(r)$. With these structures $R$ is a right
$\hH$-comodule. Similarly, $L$ is a right $\hH$-comodule via the right
regular $L$-action and $\hH_L$-coaction $L\to L\ot_L H\cong H$, $l\mapsto 1_L
\ot_L t_L(l)\cong t_L(l)$, the right $R$-action $L\ot R \to L$, $l\ot r\mapsto
\varepsilon_L\big(s_R(r)\big)l$ and $\hH_R$-coaction $L \to L\ot_R H\cong H$,
$l\mapsto 1_L \ot_R t_L(l)\cong t_L(l)$. The map $\varepsilon_R \circ t_L:L\to
R$ is an isomorphism of $\hH$-comodules with the inverse $\varepsilon_L\circ
s_R$. 

By \cite[28.4]{BrzWis:coring}, for any right $\hH$-comodule $M$ there are
isomorphisms $M^{co\hH_R}\cong \mathrm{Hom}^{\hH_R}(R,M)$ and $M^{co\hH_L}\cong
\mathrm{Hom}^{\hH_L}(L,M)$. Therefore, the following sequence of isomorphisms
and inclusions holds.
$$
M^{co\hH_R}\cong \mathrm{Hom}^{\hH_R}(R,M)\subseteq  \mathrm{Hom}^{\hH_L}(R,M)
\cong \mathrm{Hom}^{\hH_L}(L,M) \cong M^{co\hH_L}.
$$
The inclusion in the second step follows by \prref{prop:ff}.
If the antipode is bijective then the same reasoning can be applied to the
opposite Hopf algebroid to conclude that also $M^{co\hH_L}\subseteq
M^{co\hH_R}$. 
\end{proof}
}

{\color{blue}
In order to have a meaningful notion of a comodule algebra of a Hopf
algebroid $\hH$, the category of $\hH$-comodules has to be monoidal. The
following was proven in \cite[{\em Corrigendum}, Theorem
  6]{BohmBrz:hgd.cleft}.   

\begin{theorem}\label{thm:hgd.com.mon}
For any Hopf algebroid ${\mathcal H}$, ${\cM}^{\mathcal H}$ is a  monoidal
category. Moreover, there are strict monoidal forgetful functors rendering
commutative the following diagram:
$$
\xymatrix{
{\cM}^{\hH}\ar[r]\ar[d]& {\cM}^{\hH_R}\ar[d]\\
{\cM}^{\hH_L}\ar[r]& {}_R{\cM}_R \ .
}
$$
\end{theorem}

\begin{proof}
The functor on the right hand side appeared already in
\thref{thm:bgd.com_cat_mon}. 
Let us explain first what is meant by the functor in the bottom row. A right
$\hH_L$-comodule $N$ is a priori a right $L$-module, and it is made an
$L$-bimodule via the left action $ln:= n_{[0]}\varepsilon_L(n_{[1]} t_L(l))$,
for $l\in L$ and $n\in N$, cf. a symmetrical form of
\equref{eq:left_R_mod}. The functor in the bottom row takes the
$\hH_L$-comodule $N$ to the $R$-$R$ bimodule $N$, with actions
\begin{equation}\eqlabel{eq:RactionsL}
r\blacktriangleright n \blacktriangleleft r':= 
\varepsilon_L\big(s_R(r')\big) n\, \varepsilon_L\big(s_R(r)\big) \equiv
n_{[0]} \varepsilon_L\big( s_R(r) n_{[1]} s_R(r')\big),\qquad 
\textrm{for } r,r'\in R,\ n\in N.
\end{equation}
In order to see commutativity of the diagram, take a right $\hH$-comodule
$M$. Composing the functor on the left hand side with the functor in the
bottom row, it takes $M$ to an $R$-$R$ bimodule with actions in
\equref{eq:RactionsL}. Applying to $M$ the functor in the top row and the
functor on the right hand side, we obtain the $R$-$R$ bimodule $M$ with
actions
\begin{equation}\eqlabel{eq:RactionsR}
rmr'= m^{[0]} \varepsilon_R\big (s_R(r) m^{[1]}\big) r'
=  m^{[0]} \varepsilon_R\big (s_R(r) m^{[1]}s_R(r') \big)
,\qquad
\textrm{for } r,r'\in R,\ m\in M,
\end{equation}
cf. \equref{eq:left_R_mod}.
By definition, the $\hH_L$-coaction on $M$ is right $R$-linear (with respect
to the action on $M$ denoted by juxtaposition). Therefore,
\begin{eqnarray*}
(rmr')_{[0]}\stac L(rmr')_{[1]}
&=& {m^{[0]}}_{[0]} \stac L {m^{[0]}}_{[1]} s_R\left(\varepsilon_R\big (s_R(r)
  m^{[1]} s_R(r')\big)\right) \\
&=& m_{[0]} \stac L {m_{[1]}}^{[1]} s_R\left(\varepsilon_R\big (s_R(r)
  {m_{[1]}}^{[2]} s_R(r')\big)\right)
=m_{[0]} \stac L s_R(r) {m_{[1]}} s_R(r').
\end{eqnarray*}
The first equality follows by \equref{eq:RactionsR} and the right $R$-linearity
of the $\hH_L$-coaction on $M$.
The second equality follows by the right $\hH_R$-colinearity
of the $\hH_L$-coaction on $M$.
By multiplicativity, right $R$-linearity and unitality of $\Delta_R$,
$\Delta_R(s_R(r))= 1_H \ot_R s_R(r)$.
In the last equality we used this identity and multiplicativity of
$\Delta_R$, and counitality of $\Delta_R$. Applying $M \ot_L \varepsilon_L$ to
both sides, we conclude that 
$$
rmr'= m_{[0]}\varepsilon_L\big(s_R(r) {m_{[1]}} s_R(r')\big)\equiv 
r\blacktriangleright m \blacktriangleleft r'.
$$
This proves commutativity of the diagram in the theorem.

Strict monoidality of the functors on the right hand side and in the bottom
row follows by \thref{thm:bgd.com_cat_mon} (and its application to 
the opposite of the bialgebroid $\hH_L$). 
In order to see strict monoidality of the remaining two functors, recall that 
by \thref{thm:bgd.com_cat_mon} (applied to $\hH_R$ and the opposite of
$\hH_L$), the $R$-module tensor product of any two $\hH$-comodules is an
$\hH_R$-comodule and an $\hH_L$-comodule, via the diagonal coactions,
cf. \equref{eq:prod_coac}. 
Compatibility of these coactions in the sense of \deref{def:hgd_comod} is
checked as follows. For $M,N\in \cM^\hH$, $m\in M$ and $n\in N$,
\begin{equation}\eqlabel{eq:dia_comp}
(m\stac R n)^{[0]}\stac R {(m\stac R n)^{[1]}}_{(1)} \stac L {(m\stac R
    n)^{[1]}}_{(2)} =
(m^{[0]} \stac R n^{[0]}) \stac R {m^{[1]}}_{(1)} {n^{[1]}}_{(1)} \stac L
  {m^{[1]}}_{(2)} {n^{[1]}}_{(2)}. 
\end{equation}
Moreover, for any $m\in M$ and $k\in H$, there is a well defined
(i.e. $R$-balanced and $L$-balanced) map 
$$
N \stac R H \stac L H \to (M\stac R N) \stac R H \stac L H,\qquad 
n\stac R h\stac L h'\mapsto (m\stac R n)\stac R k_{(1)} h\stac L k_{(2)} h',
$$
where we used that the range of the coproduct of $\hH_L$ lies within the
Takeuchi product $H \times_L H$. 
Composing it with the equal maps $N \to N \ot_R H \ot_L H$, $n\mapsto 
n^{[0]}\ot_R {n^{[1]}}_{(1)} \ot_ L {n^{[1]}}_{(2)} =
{n_{[0]}}^{[0]} \ot_R {n_{[0]}}^{[1]} \ot_L n_{[1]}$, we conclude that the
right hand side of \equref{eq:dia_comp} is equal to 
$(m^{[0]} \ot_R {n_{[0]}}^{[0]}) \ot_R {m^{[1]}}_{(1)} {n_{[0]}}^{[1]} \ot_L
{m^{[1]}}_{(2)} n_{[1]}$.
Similarly, since the range of the $\hH_R$-coaction on $N$ lies within the
Takeuchi product $N\times_R H$, for any $n\in N$ and $k\in H$ there is a well
defined map 
$$
M \stac R H \stac L H \to (M\stac R N) \stac R H \stac L H,\qquad 
m\stac R h\stac L h'\mapsto (m\stac R n^{[0]})\stac R h n^{[1]}\stac L h' k.
$$
Composing it with the equal maps $M \to M \ot_R H \ot_L H$, $m\mapsto 
m^{[0]}\ot_R {m^{[1]}}_{(1)} \ot_ L {m^{[1]}}_{(2)} =
{m_{[0]}}^{[0]} \ot_R {m_{[0]}}^{[1]} \ot_L m_{[1]}$, we conclude that the
right hand side of \equref{eq:dia_comp} is equal also to 
$$
({m_{[0]}}^{[0]} \stac R {n_{[0]}}^{[0]}) \stac R
{m_{[0]}}^{[1]}{n_{[0]}}^{[1]} \stac L m_{[1]} n_{[1]}=
{(m\stac R n)_{[0]}}^{[0]} \stac R {(m\stac R n)_{[0]}}^{[1]} \stac L {(m\stac R
  n)_{[1]}} .
$$
The other compatibility relation in \deref{def:hgd_comod} is checked by
similar steps.
Recall from the proof of \coref{cor:hgd.coinv} that $R$($\cong L$) is a right
$\hH$-comodule as well. 
Finally, the $R$-module tensor product of $\hH$-comodule maps is an
$\hH_R$-comodule map and an $\hH_L$-comodule map by
\thref{thm:bgd.com_cat_mon}. Thus it is an $\hH$-comodule map. By Theorem 
\thref{thm:bgd.com_cat_mon} also the coherence natural transformations in
${}_R \cM_R$ are $\hH_R$- and $\hH_L$-comodule maps, so $\hH$-comodule maps,
what completes the proof.
\end{proof}

In light of Theorem \ref{thm:hgd.com.mon}, comodule algebras of a Hopf
algebroid are defined as follows.

\begin{definition}
A right {\em comodule algebra} of a Hopf algebroid $\hH$ is a monoid in the
monoidal category $\cM^\hH$ of right $\hH$-comodules. Explicitly, this means
an $R$-ring and right $\hH$-comodule $A$, such that the unit map $R \to A$ and
the multiplication map $A\ot_R A \to A$ are right 
$\hH$-comodule maps. Using the notations $a\mapsto a^{[0]}\ot_R a^{[1]}$ and
$a\mapsto a_{[0]}\ot_L a_{[1]}$ for the $\hH_R$- and $\hH_L$-coactions,
respectively, $\hH$-colinearity of the unit and the multiplication means the
following identities, for all $a,a'\in A$:
\begin{eqnarray*}
&{1_A}^{[0]} \stac  R {1_A}^{[1]} = 1_A\stac  R 1_H,\qquad 
&(aa')^{[0]} \stac  R (aa')^{[1]} = a^{[0]} a^{\prime [0]} \stac  R a^{[1]}
  a^{\prime [1]}\\
&1_{A[0]} \stac  L 1_{A{[1]}} = 1_{A}\stac  L 1_H,\qquad 
&(aa')_{[0]} \stac  L (aa')_{[1]} = a_{[0]} a'_{[0]} \stac  L a_{[1]}
  a'_{[1]} .
\end{eqnarray*}
\end{definition}

\begin{definition}
For an $\hH$-comodule algebra $A$ and $B:=A^{co\hH_R}$, we say that
$B\subseteq A$ is an $\hH$-extension.
\end{definition}
}

The related $\hH_R$-, and $\hH_L$-coactions in an $\hH$-extension 
$B\subseteq A$ determine two canonical maps,
\begin{eqnarray}
&\can_R:A\stac B A \to A\stac R H,\qquad 
&a\stac B a'\mapsto aa^{\prime[0]}\stac
R a^{\prime[1]}\qquad \textrm{and}\nonumber\\ 
&\can_L:A\stac B A \to A\stac L H,\qquad 
&a\stac B a'\mapsto a_{[0]}a'\stac L\eqlabel{eq:hgd.cans}
a_{[1]}.
\end{eqnarray}
Bijectivity of the maps in \equref{eq:hgd.cans} 
{\color{blue} 
implies that} $B\subseteq A$ is a Galois extension by the right 
bialgebroid $\hH_R$, and the left bialgebroid $\hH_L$, respectively. It is not
known in general if bijectivity of one implies bijectivity of the other. A
partial result is given in following \cite[Lemma 3.3]{Bohm:hgdGal}. 

\begin{proposition}\prlabel{prop:hgd.cans}
For an algebra extension $B\subseteq A$ by a Hopf algebroid $\hH$ with a {\em
  bijective} antipode $S$, the canonical map $\can_L$ in \equref{eq:hgd.cans}
  is bijective if and only if $\can_R$ is bijective.
 \end{proposition}
\begin{proof} (Sketch.) 
By the bijectivity of the antipode, for any right $\hH$-comodule $M$, the map
$\Phi_M$
in \equref{eq:Phi} is bijective, with the inverse $m\otimes_L h \mapsto
m^{[0]}\otimes_R S^{-1}(h)  
m^{[1]}$. The canonical maps in \equref{eq:hgd.cans} are related by
$\Phi_A\circ \can_R =\can_L$, what proves the claim.
\end{proof}

{\color{blue}
The Galois theory of a Hopf algebroid $\hH$ is greatly simplified whenever
the category of $\hH$-comodules is isomorphic to the categories of
comodules of the constituent left and right bialgebroids: In this case
$\cM^\hH$ can be described as a category of an appropriate coring.
The following was obtained in \cite[\emph{Corrigendum}, Theorem
  4]{BohmBrz:hgd.cleft}.}  

{\color{blue}
\begin{theorem}\thlabel{thm:pure}
Consider a Hopf algebroid $\hH$, with structure maps denoted as in
\deref{def:hgd}, and the forgetful functors $F_L:\cM^{\hH_L}\to \cM_k$,
$F_R:\cM^{\hH_R}\to \cM_k$, $G_L:\cM^\hH\to \cM^{\hH_L}$ and $G_R:\cM^\hH\to
\cM^{\hH_R}$.

(1) If the equalizer 
\begin{equation}\eqlabel{eq:R.M}
\xymatrix{
M \ar[rr]^-{\varrho_R}&&
M\stac R H \ar@<2pt>[rr]^-{\varrho_R\ot_R H}\ar@<-2pt>[rr]_-{M\ot_R \Delta_R}&&
M\stac R H \stac R H
}
\end{equation}
in $\cM_L$ is $H\ot_L H$-pure, i.e. it is preserved by the functor $- \ot_L
H\ot_L H:\cM_L\to \cM_L$, for any right $\hH_R$-comodule $(M,\varrho_R)$, then
there exists a functor $U:\cM^{\hH_R}\to \cM^{\hH_L}$, such that $F_L \circ U
=F_R$ and $U\circ G_R=G_L$. 

(2) If the equalizer 
\begin{equation}\eqlabel{eq:L.M}
\xymatrix{
N \ar[rr]^-{\varrho_L}&&
N\stac L H \ar@<2pt>[rr]^-{\varrho_L\ot_L H}\ar@<-2pt>[rr]_-{N\ot_L \Delta_L}&&
N\stac L H \stac L H
}
\end{equation}
in $\cM_R$ is $H\ot_R H$-pure, i.e. it is preserved by the functor $-\ot_R
H\ot_R H:\cM_R\to \cM_R$, for any right $\hH_L$-comodule $(N,\varrho_L)$, then
there exists a functor $V:\cM^{\hH_L}\to \cM^{\hH_R}$, such that $F_R \circ V
=F_L$ and $V\circ G_L=G_R$. In particular, $G_L$ is full.

(3) If both purity assumptions in parts (1) and (2) hold, then the forgetful
functors $G_R:\cM^\hH\to \cM^{\hH_R}$ and $G_L:\cM^\hH\to \cM^{\hH_L}$
are isomorphisms, hence $U$ and $V$ are inverse isomorphisms.
\end{theorem}

\begin{proof}
(1) 
Recall that \equref{eq:R.M} defines the $\hH_R$-cotensor product $M
\Box_{\hH_R} H\cong M$. By axiom (ii) in \deref{def:hgd}, $H$ is an
$\hH_R$-$\hH_L$ bicomodule, with left coaction $\Delta_R$ and right coaction
$\Delta_L$. Thus in light of \cite[22.3]{BrzWis:coring} and its Erratum, we can
define a desired functor $U:= - \Box_{\hH_R} H$. Clearly, it satisfies $F_L
\circ U = F_R$. 
For an $\hH$-comodule $(M,\varrho_L,\varrho_R)$, the coaction on the
$\hH_L$-comodule $U\big( G_R(M,\varrho_L,\varrho_R)\big)=U(M,\varrho_R)$ is
given by  
\begin{equation}\eqlabel{eq:rho'}
\xymatrix{
M \ar[r]^-{\varrho_R}&
M\Box_{\hH_R} H \ar[rr]^-{M\Box_{\hH_R} \Delta_L}&&
M\Box_{\hH_R} (H\stac L H) \ar[r]^-{\simeq}&
(M\Box_{\hH_R} H)\stac L H \ar[rr]^-{M\ot_R \varepsilon_R \ot_L H}&&
M\stac L H,
}
\end{equation}
where in the third step we used that since the equalizer \equref{eq:R.M} is
$H\ot_L H$-pure, it is in particular $H$-pure. Using that $\varrho_R$ is a
right $\hH_L$-comodule map and counitality of $\varrho_R$, we conclude that
\equref{eq:rho'} is equal to $\varrho_L$. Hence $U\circ G_R =G_L$. 
Note that this yields an alternative proof of fully faithfulness of $G_R$.
Indeed,
this proves that for any two $\hH$-comodules $M$ and $M'$, and any
$\hH_R$-comodule map $f:M\to M'$, $U(f)=f$ is an $\hH_L$-comodule map hence
an $\hH$-comodule map. Part (2) is proven symmetrically. 

(3) 
For the functor $U$ in part (1) and a right $\hH_R$-comodule $(M,\varrho_R)$,
denote $U(M,\varrho_R)=:(M,\varrho_L)$. 
With this notation, define a functor ${G}_R^{-1}:\cM^{\hH_R}\to
\cM^\hH$, with object map $(M,\varrho_R)\mapsto (M,\varrho_R,\varrho_L)$, and
acting on the morphisms as the identity map. Being coassociative, $\varrho_R$
is an $\hH_R$-comodule map, so by part (1), $U(\varrho_R)=\varrho_R$ is an
$\hH_L$-comodule map. Symmetrically, by part (2), $V(\varrho_L)=\varrho_L$ is
an $\hH_R$-comodule map. So 
$G_R^{-1}$ is a well defined functor. One easily checks that it is the inverse
of $G_R$. 

In a symmetrical way, in terms of the functor $V(N,\varrho_L)=:(N,\varrho_R)$
in part (2), one constructs $G_L^{-1}$ with object map $(N,\varrho_L)\mapsto
(N,\varrho_L,\varrho_R)$, and acting on the morphisms as the identity map. 
The identities $G_L \circ G_R^{-1}=U$ and $G_R \circ G_L^{-1}=V$ prove that
$U$ and $V$ are mutually inverse isomorphisms, as stated.
\end{proof}

\begin{definition}
Hopf algebroids for that the purity assumptions in \thref{thm:pure} (3) hold,
are termed {\em pure} Hopf algebroids.
\end{definition}

All known examples of Hopf algebroids are pure, cf. \cite[\emph{Corrigendum},
  Example 5]{BohmBrz:hgd.cleft}.
}

\begin{example}\exlabel{ex:hgdGal}
(1) \exref{ex:Hopf.Gal} (1) can be extended as follows. Consider a Hopf
  algebroid $\hH$ and use the notations introduced in and after
  \deref{def:hgd}. The coproduct $\Delta_R$ in $\hH_R$ equips the total
  algebra $H$ with a right $\hH_R$-comodule algebra structure and the
  coproduct $\Delta_L$ in $\hH_L$ equips $H$ with a right $\hH_L$-comodule
  algebra structure. In this way $H$ becomes an $\hH$-comodule algebra.
$\hH_R$-coinvariants are those elements $h\in H$, for
  which $h^{(1)}\otimes_R h^{(2)}=h\otimes_R 1_H$, i.e. the image of $R^{op}$
  under $t_R$. $\hH_L$-coinvariants are those elements $h\in H$, for
  which $h_{(1)} \ot_L h_{(2)} = h\ot_L 1_H$, i.e. elements of
  $t_R(R^{op})=s_L(L)$. 
The canonical map 
$$
H\stac {L} H \to H\stac R H,\qquad h\stac {L} h'\mapsto h
h^{\prime(1)}\stac R h^{\prime(2)}
$$
is bijective, with the inverse $h\otimes_R h'\mapsto h
S(h'_{(1)})\otimes_{L} h'_{(2)}$. That is, the algebra extension
$L\subseteq H$, given by $s_L$ (equivalently, the algebra extension
$R^{op}\subseteq H$, given by $t_R$) is a Galois extension by $\hH_R$. In
other words, a constituent right bialgebroid $\hH_R$ in a Hopf algebroid $\hH$
provides an example of a (right) {\em $\times_R$-Hopf algebra}, in the sense
of \cite{Sch:dua}.  

The other canonical map 
$$
H \stac L H \to H\stac L H, \qquad h\stac L h'\mapsto h_{(1)} h'\stac L h_{(2)}
$$
(mind the different $L$-actions in the domain and the codomain!)
is bijective provided that $S$ is bijective. In this case the inverse is given
by $h'\ot_L h \mapsto h^{(2)} \ot_L S^{-1}(h^{(1)}) h'$.

(2) Consider a right depth 2 and balanced algebra extension $B\subseteq A$, as
in \exref{ex:depth2}. We have seen in \exref{ex:depth2} that $B\subseteq A$ is
a Galois extension by a right bialgebroid $\hH_R$.
Assume that the extension $B\subseteq A$ is also {\em Frobenius} (i.e. $A$
is a finitely generated and projective right $B$-module and
$\mathrm{Hom}_B(A,B)\cong A$ as $B$-$A$ bimodules). In this situation, the
right bialgebroid $\hH$ was proven to be a constituent right bialgebroid in a
Hopf algebroid $\hH$ in \cite[Section 3]{BohmSzl:hgdax}.
{\color{blue}
The Hopf algebroid $\hH$ is finitely generated and projective as a left
$L$-module and as a left $R$-module. This implies that $\hH$ is a pure Hopf
algebroid, hence a right $\hH_R$-comodule algebra $A$ is also a right
$\hH$-comodule algebra.
Moreover, the antipode of $\hH$ is bijective, hence the $\hH_R$-Galois
extension $B\subseteq A$ is also $\hH_L$-Galois by \prref{prop:hgd.cans}.
}
\end{example}

\section{Cleft extensions}\selabel{sec:cleft}

In this section a particular class of Galois extensions, so called cleft
extensions will be studied. These are the simplest and best understood
examples of Galois extensions and also the closest ones to the classical
problem of Galois extensions of fields.

Similarly to \seref{sec:defs}, we start with reviewing the most classical
case of a cleft extension by a Hopf algebra in \cite{DoiTak:cleft} and
\cite{BlatMont:Crprod&Gal}. A cleft extension by a coalgebra was introduced
in \cite{BrzMaj:coabund} and 
\cite{Brz:coalg.cleft} and studied further in \cite{Abu:Mor.cor} and
\cite{CaVerWang:Mor.cor.cleft}. In all these papers the analysis is based on
the study of an {\em entwining structure} associated to a Galois
extension. Here we present an equivalent description of a cleft extension by 
a coalgebra, which avoids using entwining structures. Instead, we make 
use of the {\em coring extension} behind. The advantage of this approach,
developed in \cite{BohmVer:Mor&cleft}, is that it provides a uniform approach
to 
{\color{blue} many kinds of}
cleft extensions, including cleft extensions by 
{\color{blue} pure}
Hopf algebroids, where the
coring extension in question does not come from an entwining structure.  

\subsection{Cleft extensions by Hopf algebras}\selabel{sec:Hopf.cleft}

The notion of a cleft extension by a Hopf algebra emerged already in papers by
Doi and Sweedler, but it became relevant by results in \cite{DoiTak:cleft} and
\cite{BlatMont:Crprod&Gal}. Note that, for a $k$-algebra $(A,\mu,\eta)$ and a
$k$-coalgebra $(C,\Delta,\epsilon)$, the $k$-module $\mathrm{Hom}_k(C,A)$ is a
$k$-algebra via the {\em convolution product}
$$
(f,g)\mapsto \mu\circ (f\stac k g)\circ \Delta,\qquad \textrm{for }f,g\in
\mathrm{Hom}_k(C,A), 
$$
and unit element $\eta\circ \epsilon$.

\begin{definition}\delabel{def:Hopf.cleft}
An algebra extension $B\subseteq A$ by a Hopf algebra $H$ is said to be {\em
  cleft} provided that there exists a convolution invertible right
  $H$-comodule map $j:H\to A$, called a {\em cleaving map}.
\end{definition}

By antipode axioms \equref{eq:Hopf.antip.ax} in a Hopf algebra $H$, the
antipode is convolution inverse of the ($H$-colinear) identity map $H$. Thus
any $k$-Hopf algebra $H$ is an $H$-cleft extension of $k$ (via the unit map).

Following \cite[Theorem 9]{DoiTak:cleft} explains in what sense cleft
extensions are distinguished Hopf Galois extensions.

\begin{theorem}\thlabel{thm:cleft=Gal&norm.b}
An algebra extension $B\subseteq A$ by a Hopf algebra $H$ is cleft if and only
if it is an $H$-Galois extension and the {\em normal basis property} holds,
i.e. $A\cong B\otimes_k H$ as left $B$-modules right $H$-comodules.
\end{theorem}

The construction of a crossed product with a bialgebra $H$ (with coproduct
$\Delta(h)=h_{(1)}\otimes_k h_{(2)}$ and counit $\epsilon$) was introduced in
\cite{BlatCohMont:Cr.prod}, as follows. 

\begin{definition} \delabel{def:Hopf.cr.prod}
A $k$-bialgebra $H$ {\em measures} a $k$-algebra $B$ if there exists a
$k$-module map $\cdot: H\otimes_k B\to B$, such that $h\cdot 1_B=\epsilon(h)
1_B$ and $h\cdot (bb')=(h_{(1)}\cdot b)(h_{(2)}\cdot b')$, for $h\in H$,
$b,b'\in B$. 

A $B$-valued {\em 2-cocycle} on $H$ is a $k$-module map $\sigma: H\otimes_k
H\to B$, such that $\sigma(1_H,h)=\epsilon(1_H)1_B =\sigma(h,1_H)$, for $h\in
H$, and
$$
\big(h_{(1)}\cdot\sigma(k_{(1)},m_{(1)})\big)\sigma(h_{(2)}, k_{(2)} m_{(2)})
= \sigma(h_{(1)},k_{(1)})\sigma(h_{(2)}k_{(2)}, m),\qquad \textrm{for }
h,k,m\in H. 
$$
The $H$-measured algebra $B$ is a {\em $\sigma$-twisted $H$-module} if in
addition
$$
\big(h_{(1)}\cdot(k_{(1)}\cdot b)\big)\sigma(h_{(2)},k_{(2)})=
\sigma(h_{(1)},k_{(1)}) (h_{(2)}k_{(2)}\cdot b),\qquad \textrm{for }
b\in B,\ h,k\in H. 
$$
\end{definition}

\begin{proposition}\prlabel{prop:Cr.prod}
Consider a $k$-bialgebra $H$ and an $H$-measured $k$-algebra $B$.
Let $\sigma:H\otimes_k H\to B$ be a $k$-module map.
The $k$-module $B\otimes_k H$ is an algebra, with multiplication 
$$
(b\#h)(b'\#h')=b(h_{(1)}\cdot b')\sigma(h_{(2)},h'_{(1)})\#
h_{(3)}h'_{(2)}
$$
and unit element $1_B\# 1_H$, if and only if $\sigma$ is a $B$-valued
2-cocycle on $H$ and $B$ is a $\sigma$-twisted $H$-module.
This algebra is called the {\em crossed
  product} of $B$ with $H$, with respect to the cocycle $\sigma$. It is
denoted by $B\#_\sigma H$.
\end{proposition}

It is straightforward to see that a crossed product algebra $B\#_\sigma H$ is 
a right $H$-comodule algebra, with coaction given in terms of the coproduct
$\Delta$ in $H$, as $B\otimes_k \Delta$. What is more, $B\subseteq B\#_\sigma
H$ is an extension by $H$. It is most natural to ask what $H$-extensions arise
as crossed products. One implication in forthcoming \thref{thm:cleft=cr.prod}
was proven first in \cite[Theorem 11]{DoiTak:cleft}. Other implication was
proven in \cite[Theorem 1.18]{BlatMont:Crprod&Gal}.
Since for a $k$-bialgebra $H$ also $H\otimes_k H$ is a $k$-bialgebra,
convolution invertibility of a $B$-valued 2-cocycle $\sigma$ on $H$ is
understood in the convolution algebra $\mathrm{Hom}_k(H\otimes_k H, B)$. 

\begin{theorem}\thlabel{thm:cleft=cr.prod}
An algebra extension $B\subseteq A$ by a Hopf algebra $H$ is a cleft extension
if and only if $A$ is isomorphic to $B\#_\sigma H$, as a left $B$-module and
right $H$-comodule algebra, for some convolution invertible $B$-valued
2-cocycle $\sigma$ on $H$. 
\end{theorem}

Another important aspect of cleft extensions is that they provide examples of
Hopf Galois extensions $B\subseteq A$, beyond the case when $A$ is a
faithfully flat $B$-module, when a Strong Structure Theorem holds. Recall that,
for a $k$-bialgebra $H$ and its right comodule algebra $A$, right-right
relative Hopf modules are right modules for the monoid $A$ in $\cM^H$. That
is, right $A$-modules and right $H$-comodules $M$, such that the $A$-action
$M\otimes_k A \to M$ is a right $H$-comodule map with respect to the diagonal
$H$-coaction in $M\otimes_k A$. Equivalently, the $H$-coaction $M\to
M\otimes _k H$ is a right $A$-module map with respect to the right $A$-action
in $M\otimes_k H$ given by the $H$-coaction in $A$. Clearly, for any right
$B:=A^{coH}$-module $N$, the tensor product $N\otimes_B A$ inherits a relative
Hopf module structure of $A$. 

\begin{theorem}\thlabel{thm:Hopf.cleft.str.str.thm}
For a cleft extension $B\subseteq A$ by a Hopf algebra $H$, the category of
right $B$-modules is equivalent to the category of right-right
$(H,A)$-relative Hopf modules, via the induction functor $-\otimes_B
A:\cM_B\to \cM^H_A$. 
\end{theorem}

\subsection{Cleft bicomodules for 
\underline{pure}
coring extensions}\selabel{sec:cleft.bicom} 

The aim of the current section is to reformulate \deref{def:Hopf.cleft} of a
cleft extension by a Hopf algebra, using the Morita theory of 
{\color{blue} pure}
coring extensions
developed in \cite{BohmVer:Mor&cleft}
{\color{blue} (see the corrected versions)}. 
This will allow us to place the results in \thref{thm:cleft=Gal&norm.b} and
\thref{thm:Hopf.cleft.str.str.thm} into a broader context. More importantly,
it will provide us with a tool of generalizations in later sections.

Observe first that, for a right comodule algebra $A$ of a $k$-bialgebra $H$,
the tensor product $\cC:=A\otimes_k H$ is an $A$-coring. The $A$-$A$ bimodule
structure is given in terms of the right $H$-coaction in $A$,
$\varrho^A:a\mapsto 
a_{[0]}\otimes_k a_{[1]}$, as 
$$
a_1(a\stac k h)a_2:= a_1 a {a_2}_{[0]}\stac k {a_2}_{[1]}, \qquad \textrm{for
} a_1,a_2\in A,\ a\stac k h\in A\stac k H.
$$
The coproduct is $A\otimes_k \Delta:A \otimes_k H\to A\otimes_k H \otimes_k
H\cong (A\otimes_k H)\otimes_A (A\otimes_k H)$, determined by the coproduct
$\Delta$ in $H$, and the counit is $A\otimes_k \epsilon:A\otimes_k H\to A$,
coming from the counit $\epsilon$ in $H$. What is more, $\cC$ is a right
$H$-comodule via 
the coaction $A\otimes_k \Delta:\cC\to \cC\otimes_k H$. The coproduct in $\cC$ 
is a right $H$-comodule map, that is, (the constituent coalgebra in) $H$ is a
right extension of $\cC$. 
{\color{blue}
Moreover, since for any right $\cC$-comodule $M$ the equalizer
\equref{eq:pure_eq} in $\cM_k$ is split by the map $M\ot_k H\ot_k \varepsilon:
M\ot_A \cC\ot_A \cC\cong M\ot_k H\ot_k H\to M\ot_A \cC \cong M\ot_k H$, 
it is preserved by any functor of domain $\cM_k$. Therefore, $H$ is pure right
extension of $\cC$.}
Via the right regular $A$-action and coaction
$\varrho^A:A\to A\otimes_k H\cong A\otimes_A \cC$, $A$ is a right
$\cC$-comodule. 

By the above motivation, we turn to a study of any $L$-coring $\cD$, 
which is a 
{\color{blue} pure}
right extension of an $A$-coring $\cC$, and an $L$-$\cC$
bicomodule $\Sigma$. 
In \cite[Proposition 3.1]{BohmVer:Mor&cleft} to any such bicomodule $\Sigma$
a Morita context $\mathbb{M}(\Sigma)$ was associated. In order to write it up
explicitly, recall that by 
{\color{blue} the purity assumption}, 
$\Sigma\cong \Sigma\Box_\cC \cC$ is also a right $\cD$-comodule. Put
$T:=\mathrm{End}^\cC(\Sigma)$. It is an $L$-ring. 
Introduce the following index notations. For the coproduct in $\cC$, write
$c\mapsto c^{(1)}\otimes_A c^{(2)}$. For the coproduct in $\cD$, write
$d\mapsto d_{(1)}\otimes_L d_{(2)}$. For the $\cC$-coaction in $\Sigma$ write
$x\mapsto x^{[0]}\otimes_A x^{[1]}$ and for the corresponding $\cD$-coaction
in $\Sigma$ write $x\mapsto x_{[0]}\otimes_L x_{[1]}$. In each case implicit
summation is understood. Then 
\begin{equation}\eqlabel{eq:Mor}
\mathbb{M}(\Sigma)=\big(
{}_L\mathrm{Hom}_L(\cD,T)\ ,\ 
{}^\cC\mathrm{End}^\cD(\cC)^{op}\ ,\ 
{}_L\mathrm{Hom}^\cD(\cD,\Sigma)\ ,\ 
\widetilde{Q}\ ,\ 
\blacklozenge\ ,\  
\lozenge\big),
\end{equation}
where 
$$
\widetilde{Q}\colon=\{\ q\in {}_A{\rm Hom}_L(\cC,\mathrm{Hom}_A(\Sigma,A))\ |\ 
c^{(1)}q(c^{(2)})(x)=q(c)(x^{^{[0]}})x^{[1]},\quad \forall x\in \Sigma, c\in
\cC\ \}.  
$$
The algebra structures, bimodule structures and connecting
homomorphisms are given by the following formulae.
\begin{eqnarray*}
&&(vv')(d)=v(d_{(1)}) v'(d_{(2)})\\
&&(uu')(c)=u' \big(u(c)\big)\\
&&(vp)(d)=v(d_{(1)})\big(p(d_{(2)})\big)\\
&&(pu)(d)=p(d)^{[0]}\epsilon_\cC\big(u(p(d)^{[1]})\big)\\
&&(qv)(c)=q(c_{[0]})v(c_{[1]})\\
&&(uq)(c)=q\big(u(c)\big)\\
&&(q \smallblackdiamond p)(c)=
c^{(1)}q({c^{(2)}}_{[0]})\big(p({c^{(2)}}_{[1]})\big)\equiv 
q(c_{[0]})\big(p({c_{[1]}})^{[0]}\big)p({c_{[1]}})^{[1]}\\  
&&(p\smalldiamond q)(d)=p(d)^{[0]}q\big(p(d)^{[1]}\big)(-),
\end{eqnarray*}
for $v,v'\in {}_L{\rm Hom}_L(\cD,T)$, $u,u'\in {}^\cC{\rm
End}^\cD(\cC)$, $p\in {}_L{\rm Hom}^\cD(\cD,\Sigma)$, $q\in \widetilde{Q}$,
$d\in \cD$ and $c\in \cC$.

As it is explained in \cite[Proposition 3.1]{BohmVer:Mor&cleft}, inspite of
its involved form, the Morita context $\mathbb{M}(\Sigma)$ has a very simple
origin. For any pure coring extension $\cD$ of $\cC$, there is a functor
$U:=-\Box_\cC \cC:\cM^\cC \to \cM^\cD$. An $L$-$\cC$ bicomodule $\Sigma$
determines another functor $V:=\mathrm{Hom}^\cC(\Sigma,-)\ot_L \cD:\cM^\cC \to
\cM^\cD$. In terms of natural transformations $\mathrm{Nat}(-,-)$ between
these functors, $\mathbb{M}(\Sigma)$ is isomorphic to the Morita context 
\begin{equation}\eqlabel{eq:cat_morita}
\big(\ \mathrm{Nat}(V,V)\ ,\ \mathrm{Nat}(U,U)\ ,\ \mathrm{Nat}(V,U)\ ,\ 
\mathrm{Nat}(U,V)\ ,\ \bullet\ ,\ \circ\ \big),  
\end{equation}
where all algebra and bimodule structures and also the connecting maps are
given by opposite composition of natural transformations.

Let us compute the Morita context $\mathbb{M}(A)$ in our motivating example,
coming from an algebra extension $B\subseteq A$ by a $k$-Hopf algebra $H$. We
claim that it is isomorphic to a sub-Morita context of a (degenerate) Morita
context, in which 
both algebras and both bimodules are equal (as $k$-modules) to
$\mathrm{Hom}_k(H,A)$, and all algebra structures, bimodule structures and
connecting homomorphisms are given by the convolution product. Indeed, in the
current 
case the $k$-algebra $L$ reduces to $k$ and the coring $\cD$ reduces to the
$k$-coalgebra in $H$. The role of the bicomodule $\Sigma$ is played by $A$ and
the endomorphism algebra $T$ is isomorphic to the coinvariant subalgebra
$B=A^{coH}$. The $A$-coring $\cC$ is equal to $A\otimes_k H$. 
The isomorphism ${}_A\mathrm{Hom} (A\otimes_k H,A)\cong \mathrm{Hom}_k (H,A)$
induces an isomorphism 
\begin{eqnarray*}
\widetilde{Q}&\cong& \{\ \widetilde{q}\in \mathrm{Hom}_k (H,A)\ |\ 
\widetilde{q}(h_{(2)})_{[0]}\stac k h_{(1)} \widetilde{q}(h_{(2)})_{[1]}=
\widetilde{q}(h)\stac k 1_H,\quad \forall h\in 
H\ \}\\
&\equiv&\{\ \widetilde{q}\in \mathrm{Hom}_k (H,A)\ |\ 
\widetilde{q}(h)_{[0]}\stac k \widetilde{q}(h)_{[1]}=
\widetilde{q}(h_{(2)})\stac k S(h_{(1)}),\quad \forall h\in H\ \}.  
\end{eqnarray*}
Thus we conclude that
$\widetilde{Q}$ is isomorphic to $\mathrm{Hom}^H(H^{tw},A)$, where $H^{tw}$ is 
the $k$-module $H$, considered to be a right $H$-comodule via the twisted
coaction $h\mapsto h_{(2)}\otimes_k S(h_{(1)})$. 

The monomorphism
$$
{}^\cC\mathrm{End}^H(A\stac k H)\hookrightarrow {}_A\mathrm{End}^H(A\stac k H)
\cong \mathrm{Hom}_k (H,A),
\qquad u\mapsto {\widetilde u}:=(A\stac k \epsilon)\circ u\circ (1_A \stac k -)
$$
establishes an isomorphism (with inverse ${\widetilde u}\mapsto (\ a\otimes_k
h\mapsto a{\widetilde u}(h_{(1)})\otimes_k h_{(2)}\ )$) from
${}^\cC\mathrm{End}^H(A\stac k H)$ to
\begin{eqnarray*}
X&:=&\{\ {\widetilde u}\in
\mathrm{Hom}_k(H,A)\ |\ {\widetilde u}(h_{(2)})_{[0]} \stac k
h_{(1)}{\widetilde u}(h_{(2)})_{[1]}={\widetilde u}(h_{(1)})\stac k 
h_{(2)},\quad \forall h\in H\ \}\\
&\equiv&\{\ {\widetilde u}\in
\mathrm{Hom}_k(H,A)\ |\ {\widetilde u}(h)_{[0]}\stac k {\widetilde
  u}(h)_{[1]}={\widetilde u}(h_{(2)})\stac k S(h_{(1)})h_{(3)}, \quad \forall
h\in H \ \}.
\end{eqnarray*}
So $\mathbb{M}(A)$ is isomorphic to the Morita context
\begin{equation}\eqlabel{eq:Hopf.Mor}
\big(
\mathrm{Hom}_k(H,B)\ ,\ 
X\ ,\ 
\mathrm{Hom}^H(H,A)\ ,\ 
\mathrm{Hom}^H(H^{tw},A)\ ,\ 
\blacklozenge'\ ,\  
\lozenge'\big).
\end{equation}
A monomorphism of Morita contexts is given by the 
obvious inclusions 
$\mathrm{Hom}_k (H,B)\hookrightarrow $\break 
$\mathrm{Hom}_k (H,A)$, 
$X\hookrightarrow \mathrm{Hom}_k (H,A)$,
$\mathrm{Hom}^H (H,A) \hookrightarrow \mathrm{Hom}_k (H,A)$, 
$\mathrm{Hom}^H(H^{tw},A)\hookrightarrow \mathrm{Hom}_k (H,A)$.
Since it is well known (cf. \cite[Lemma 3.2]{Doi:cleft}) that the convolution
inverse of a right $H$-comodule map $j:H\to A$, if it exists, belongs to
$\mathrm{Hom}^H(H^{tw},A)$, we obtained the following reformulation of
\deref{def:Hopf.cleft}.
\begin{proposition}\prlabel{prop:Hopf.cleft<Mor}
An algebra extension $B\subseteq A$ by a Hopf algebra $H$ is cleft if and
only if there exist elements $j\in \mathrm{Hom}^H(H,A)$ and ${\widetilde j}\in
\mathrm{Hom}^H(H^{tw},A)$ in the two bimodules in the Morita context
\equref{eq:Hopf.Mor} which are mapped by the connecting maps $\blacklozenge'$
and $\lozenge'$ to the unit elements of the respective algebras in the Morita
context. 
\end{proposition}

Motivated by \prref{prop:Hopf.cleft<Mor} we impose following
\cite[Definition 5.1]{BohmVer:Mor&cleft}.

\begin{definition}\delabel{def:cleft.bicom}
Let an $L$-coring $\cD$ be a 
{\color{blue} pure}
right extension of an $A$-coring $\cC$. An
$L$-$\cC$ bicomodule $\Sigma$ is said to be {\em cleft} provided that there
exist elements $j\in {}_L\mathrm{Hom}^\cD(\cD,\Sigma)$ and $\widetilde{j}\in 
\widetilde{Q}$ in the two bimodules in the associated Morita context
$\mathbb{M}(\Sigma)$ in \equref{eq:Mor}, such that $j \smalldiamond
\widetilde{j}$ is equal to the 
unit element in the convolution algebra ${}_L \mathrm{Hom}_L (\cD,T)$ and
$\widetilde{j}\smallblackdiamond j$ is equal to the unit element in the other
algebra ${}^\cC\mathrm{End}^\cD(\cC)^{op}$ in the Morita context
$\mathbb{M}(\Sigma)$. We say that $j$ and $\widetilde{j}$ are {\em mutual
  inverses} in $\mathbb{M}(\Sigma)$. 
\end{definition}

In the forthcoming sections we will define cleft extensions $B\subseteq A$ by
coalgebras, corings and 
{\color{blue} pure}
Hopf algebroids, by finding an appropriate coring
extension and requiring $A$ to be a cleft bicomodule for it in the sense of
\deref{def:cleft.bicom}.

In the rest of this section we recall some results from
\cite{BohmVer:Mor&cleft} about cleft bicomodules, extending
\thref{thm:cleft=Gal&norm.b} and \thref{thm:Hopf.cleft.str.str.thm}. 

\begin{theorem}\thlabel{thm:Mor_cleft=Gal&norm.b}
Consider an $L$-coring $\cD$ which is a 
{\color{blue} pure}
right extension of an $A$-coring $\cC$.
For an $L$-$\cC$ bicomodule $\Sigma$, put $T:=\mathrm{End}^\cC(\Sigma)$.
The bicomodule $\Sigma$ is cleft if and only if the following hold.

(1) 
The natural transformation of functors $\cM_A\to \cM^\cC$
\begin{equation}\eqlabel{eq:Wisbauer.Gal}
\mathrm{can}:\mathrm{Hom}_A(\Sigma,-)\stac T
\Sigma\to -\stac A \cC,\qquad
\mathrm{can}_N:\phi_N\stac T x \mapsto \phi_N(x^{[0]})\stac A x^{[1]}
\end{equation}
is a natural isomorphism.

(2) The {\em normal basis property} holds,
i.e. $\Sigma\cong T\otimes_L \cD$ as left $T$-modules right $\cD$-comodules.
\end{theorem}

\begin{proof} (Sketch.)
Let us use the index notations introduced in the paragraph preceding
\equref{eq:Mor}. If $\Sigma$ is a cleft bicomodule then the inverse of the
natural transformation \equref{eq:Wisbauer.Gal} can be constructed in terms of
the mutually inverse elements $j$ and $\widetilde j$ in the two bimodules of
the Morita context \equref{eq:Mor}:
$$
\mathrm{can}_N^{-1}:N \stac A \cC \to \mathrm{Hom}_A(\Sigma,N)\stac T \Sigma,
\qquad n\stac A c \mapsto n {\widetilde j}(c_{[0]})(-)\stac T j(c_{[1]}).
$$
This proves property (1). In order to verify property (2),
a left $T$-module right $\cD$-comodule isomorphism $\kappa:\Sigma\to
T\otimes_L \cD$ is constructed as
$
x\mapsto {x_{[0]}}^{[0]}{\widetilde j}({x_{[0]}}^{[1]})(-)\otimes_L x_{[1]}
$,
with the inverse 
$t\otimes_L d \mapsto t(j(d))$.

Conversely, if \equref{eq:Wisbauer.Gal} is a natural isomorphism and
$\kappa:\Sigma\to T\otimes_L \cD$ is an isomorphism of left $T$-modules and
right $\cD$-comodules, then the required mutually inverse elements in the two
bimodules of the Morita context \equref{eq:Mor} are given by 
$$
j:= \kappa^{-1}(1_T\stac L -)\qquad \textrm{and}\qquad  
{\widetilde j}:=[\mathrm{Hom}_A(\Sigma,A)\stac T  (T\stac L \epsilon_\cD)\circ 
\kappa]\circ \mathrm{can}_A^{-1},
$$ 
respectively, where $\epsilon_\cD$ is the counit in the coring $\cD$.

For more details we refer to \cite[Theorem 3.6]{BohmVer:Mor&cleft}.
\end{proof}

A right $\cC$-comodule $\Sigma$, for which \equref{eq:Wisbauer.Gal} is a
natural isomorphism, was termed a {\em Galois comodule} in
\cite{Wis:Gal.com}. Note that if $\Sigma$ is a finitely generated and
projective right $A$-module then \equref{eq:Wisbauer.Gal} is a
natural isomorphism if and only if $\mathrm{can}_A$ is bijective. Hence a
right comodule algebra $A$ of a Hopf algebra $H$ is a Galois comodule for the
canonical $A$-coring $A\otimes_k H$ if and only if $A^{coH}\subseteq A$ is an
$H$-Galois extension. Thus \thref{thm:Mor_cleft=Gal&norm.b} extends
\thref{thm:cleft=Gal&norm.b}. 

\begin{theorem}\thlabel{thm:corext.cleft.str.str.thm}
Consider an $L$-coring $\cD$ which is a 
{\color{blue} pure}
right extension of an $A$-coring $\cC$.
Take a cleft $L$-$\cC$ bicomodule $\Sigma$ and put
$T:=\mathrm{End}^\cC(\Sigma)$. If there exist finite sets, $\{v_i\}\subseteq
{}_L\mathrm{Hom}_L (\cD,T)$ and $\{d_i\}\subseteq \cD$ satisfying $\sum_i
v_i(d_i)=1_T$, then the category of
right $T$-modules is equivalent to the category of right
$\cC$-comodules, via the induction functor $-\otimes_T \Sigma:\cM_T\to
\cM^\cC$.  
\end{theorem}

\begin{proof} (Sketch.)
The induction functor $-\otimes_T \Sigma:\cM_T\to \cM^\cC$ possesses a right
adjoint, the functor $\mathrm{Hom}^\cC(\Sigma,-)$. The counit of the
adjunction is given by evaluation,
$$
n_M:\mathrm{Hom}^\cC(\Sigma,M)\stac T \Sigma \to M,\qquad  \varphi_M\stac T x
\mapsto \varphi_M(x),
$$ 
for any right $\cC$-comodule $M$. Denote the right $\cC$-coaction in $M$ by
$m\mapsto m^{[0]}\otimes_A m^{[1]}$ and the corresponding 
(via the isomorphism $M\cong M\Box_\cC\,\cC$)
$\cD$-coaction by $m\mapsto m_{[0]}\otimes_L m_{[1]}$. If
$\Sigma$ is a cleft bicomodule then the inverse of the counit can be
constructed in terms of the mutually inverse elements $j$ and $\widetilde j$
in the two bimodules of the Morita context \equref{eq:Mor}:
$$
n_M^{-1}(m):= {m_{[0]}}^{[0]}{\widetilde j}({m_{[0]}}^{[1]})(-) \stac T
j(m_{[1]})\equiv\mathrm{can}_M^{-1}(m^{[0]}\stac A m^{[1]}),
$$
where in the last expression the natural isomorphism \equref{eq:Wisbauer.Gal}
appears (cf. \thref{thm:Mor_cleft=Gal&norm.b}).

The unit of the adjunction is
$$
u_N:N\to \mathrm{Hom}^\cC(\Sigma,N\stac T \Sigma),\qquad 
n\mapsto \big(\ x\mapsto n\stac T x\ \big),
$$
for any right $T$-module $N$. Any 
element $q$ of the bimodule $\widetilde Q$ in the Morita context
\equref{eq:Mor} determines 
a map $\varpi_q:\Sigma \to T
$, 
$
\varpi_{q}(x)=x^{[0]}q(x^{[1]})(-).
$
In terms of the mutually inverse elements $j$ and $\widetilde j$
in the two bimodules of the Morita context \equref{eq:Mor}, and the finite sets
$\{v_i\}\subseteq {}_L\mathrm{Hom}_L (\cD,T)$ and $\{d_i\}\subseteq \cD$ in
the Theorem, the inverse of the unit is constructed as
$$
u_N^{-1}(\zeta_N):= \sum_i[(N\stac T \varpi_{\tilde j})\circ \zeta_N] 
\big((v_i(d_{i(1)})(j(d_{i(2)}))\big).
$$
For more details we refer to Theorem 4.1, Theorem 2.6 and Proposition 4.3 in
\cite{BohmVer:Mor&cleft}. 
\end{proof}

If $B\subseteq A$ is a cleft extension by a Hopf algebra $H$, then there exist
one-element sets satisfying the condition in
\thref{thm:corext.cleft.str.str.thm}. Indeed, for $v$ one can choose the
counit in $H$ and for $d$ the unit element in $H$. Thus
\thref{thm:corext.cleft.str.str.thm} extends
\thref{thm:Hopf.cleft.str.str.thm}.

\subsection{Cleft extensions by coalgebras}\selabel{sec:coalg.cleft}
 
Consider a coalgebra $(C,\Delta,\epsilon)$ over a commutative ring $k$. 
The cleft property of a $C$-extension $B\subseteq A$ 
appeared first in a more restricted form in \cite{BrzMaj:coabund} (and was
related to a crossed product with a coalgebra in \cite{Brz:crosprod}). The
definition we use here was introduced in
\cite{Brz:coalg.cleft}, and studied further in \cite{Abu:Mor.cor} and
\cite{CaVerWang:Mor.cor.cleft}. The aim of this section is to
reformulate the definition of the cleft property of a $C$-extension
$B\subseteq A$ in the spirit of \deref{def:cleft.bicom}. That is, we
want to see that a $C$-extension $B\subseteq A$ is cleft if and only if $A$ is
a cleft bicomodule for an appropriate 
{\color{blue} pure}
coring extension. With the experience of
Hopf algebra cleft extensions in mind, the obvious candidate is the coalgebra
$C$ to be a
{\color{blue} pure} 
right extension of an $A$-coring $A\otimes_k C$. Note however,
that $A\otimes_k C$ is not an $A$-coring without further assumptions. By
\cite[Proposition 2.2]{Brz:str}, $\cC:=A\otimes_k C$ is an $A$-coring (with the
left regular $A$-module of the first factor, coproduct $A\otimes_k \Delta$ and
counit $A\otimes_k \epsilon$) if and only if the algebra $A$ and the
coalgebra $C$ are {\em entwined}. Furthermore, in this case, 
$A$ is an entwined module (i.e. a right $\cC$-comodule) if and only if 
the (given) $C$-coaction on $A$ is a right $A$-module map $A\to \cC$.
In light of \cite[Proposition 2.3 (1) $\Leftrightarrow$ (2)]{Brz:coalg.cleft},
the following is an equivalent formulation of the definition of a cleft
coalgebra extension in \cite[p. 293]{Brz:coalg.cleft}.
\begin{definition}\delabel{def:coalg.cleft}
Let $(C,\Delta,\epsilon)$ be a $k$-coalgebra and $B\subseteq A$ be an algebra
extension by $C$. The $C$-extension $B\subseteq A$ is said to be {\em
$C$-cleft} provided that the following properties hold.

(1) $\cC:=A\otimes_k C$ is an $A$-coring, with coproduct $A\otimes_k \Delta$
and counit $A\otimes_k \epsilon$, the left regular $A$-module structure of
the first factor and a right $A$-module structure such that the $C$-coaction
in $A$ is a right $A$-module map $A\to \cC$.

(2) There exists a convolution invertible right $C$-comodule map $j:C\to A$,
the so called {\em cleaving map}. 
\end{definition}

Following \cite[Proposition 6.1]{BohmVer:Mor&cleft} reveals the relation of
\deref{def:coalg.cleft} to \deref{def:cleft.bicom}.

\begin{proposition}\prlabel{pr:coalg.cleft}
Let $(C,\Delta,\epsilon)$ be a $k$-coalgebra and $B\subseteq A$ an extension
of algebras. The algebra extension $B\subseteq A$ is $C$-cleft if and only if
the following conditions hold.

(1) $\cC:=A\otimes_k C$ is an $A$-coring (with the left regular $A$-module of
the first factor, coproduct $A\otimes_k \Delta$ and counit $A\otimes_k
\epsilon$), hence $C$ is a 
{\color{blue} pure}
right extension of $\cC$.

(2) The right regular $A$-module extends to a cleft bicomodule for the coring
    extension $C$ of $\cC$.

(3) $B=A^{coC}$.
\end{proposition}

\begin{proof}(Sketch.)
$A$ is a right $\cC$-comodule if and only if it is a right $C$-comodule such
  that the coaction is a right $A$-module map $A\to\cC$. By definition,
  $B=A^{coC}$ if and only if $B\subseteq A$ is an extension by $C$.  
Analogously to \equref{eq:Hopf.Mor}, the Morita
context $\mathbb{M}(A)$, corresponding to the $k$-$\cC$ bicomodule $A$ via 
\equref{eq:Mor}, is isomorphic to
$$
\big(\ 
\mathrm{Hom}_k(C,B)\ ,\ 
X\ ,\ 
\mathrm{Hom}^C(C,A)\ ,\ 
{\widetilde Q}'\ ,\ 
\blacklozenge'\ ,\  
\lozenge'\big),
$$
where
\begin{eqnarray*}
{\widetilde Q}'&=&\{\ {\widetilde q}\in \mathrm{Hom}_k(C,A)\ |\  
(1_A\stac k c_{(1)}) {\widetilde q}(c_{(2)})= {\widetilde q}(c) 1_{A[0]}
\stac k 1_{A[1]},\quad \forall c\in C\ \}
\qquad \textrm{and}\\ 
X'&=&\{\ {\widetilde u}\in \mathrm{Hom}_k(C,A)\ |\ 
(1_A\stac k c_{(1)}) {\widetilde u}(c_{(2)})= {\widetilde
  u}(c_{(1)})\stac k c_{(2)},\quad \forall c\in C\ \},
\end{eqnarray*}
and all algebra and bimodule structures and also the connecting maps
$\blacklozenge'$ and $\lozenge'$ are given by the convolution product. Thus
the claim follows by \cite[Lemma 4.7 1.]{Abu:Mor.cor}, stating that the
convolution inverse of $j\in \mathrm{Hom}^C(C,A)$, if it exists, belongs to
${\widetilde Q}'$. 
\end{proof}

By application of \thref{thm:Mor_cleft=Gal&norm.b}, we recover a
characterization of $C$-cleft extensions in the last paragraph of Section 4 in 
\cite{Brz:coalg.cleft}, \cite[Theorem 4.9 1.$\Leftrightarrow$3.]{Abu:Mor.cor}
or \cite[Theorem 4.5 1)$\Leftrightarrow$3)]{CaVerWang:Mor.cor.cleft}.

\begin{theorem}
An algebra extension $B\subseteq A$ by a $k$-coalgebra $C$ is cleft if and
only if it is a $C$-Galois extension and the {\em normal basis property}
holds, i.e. $A\cong B\otimes_k C$ as left $B$-modules right $C$-comodules.
\end{theorem}

Application of \thref{thm:corext.cleft.str.str.thm} yields the following.

\begin{theorem}\thlabel{thm:coalg.cleft.str.str.thm}
Consider a $k$-coalgebra $C$ and a $C$-cleft algebra extension $B\subseteq
A$. 
If there exist finite sets, $\{v_i\}\subseteq \mathrm{Hom}_k(C,B)$ and 
$\{d_i\}\subseteq C$ satisfying $\sum_i v_i(d_i)=1_B$, then the category of
right $B$-modules is equivalent to the category of right comodules for the
$A$-coring $\cC:=A\otimes_k C$, via the induction functor $-\otimes_B
A:\cM_B\to \cM^\cC$.  
\end{theorem}

Consider a $k$-coalgebra $(C,\Delta,\epsilon)$ and an algebra extension
$B\subseteq A$ by $C$. If there exists a grouplike element in $C$ then it is
mapped by the $k$-module map $\epsilon(-)1_B:C\to B$ to $1_B$. Hence
\thref{thm:coalg.cleft.str.str.thm} extends \cite[Theorem 4.9
1.$\Rightarrow$2.]{Abu:Mor.cor} or \cite[Theorem 4.5
1)$\Rightarrow$2)]{CaVerWang:Mor.cor.cleft}. However, it goes beyond the
quoted theorems. The premises of \thref{thm:coalg.cleft.str.str.thm} clearly
hold whenever the counit of $C$ is surjective, e.g. if $C$ is a faithfully
flat $k$-module (cf. \cite[Theorem 4.9 1.$\Rightarrow$5.]{Abu:Mor.cor}).

\subsection{Cleft extensions by corings}\selabel{sec:coring.cleft}

Motivated by \prref{prop:Hopf.cleft<Mor} and \prref{pr:coalg.cleft}, we
propose the following definition of a cleft extension by a coring.

\begin{definition}
Let $(\cD,\Delta,\epsilon)$ be an $R$-coring and $B\subseteq A$ an extension
of algebras. The algebra extension $B\subseteq A$ is said to be {\em
  $\cD$-cleft} provided that the following conditions hold. 

(1) $\cC:=A\otimes_R \cD$ is an $A$-coring (with the left regular $A$-module
    structure of the first factor, coproduct $A\otimes_R \Delta$ and counit
    $A\otimes_R \epsilon$), hence the coring $\cD$ is a 
{\color{blue} pure}
right extension of
    $\cC$. 

(2) The right regular $A$-module extends to a cleft bicomodule for the coring
    extension $\cD$ of $\cC$.

(3) $B=A^{co\cD}$.
\end{definition}

Note that, for an $R$-ring $(A,\mu,\eta)$ and an
$R$-coring $(\cD,\Delta,\epsilon)$, the $k$-module ${}_R\mathrm{Hom}_R(\cD,A)$
is a $k$-algebra via the {\em convolution product}
$$
(f,g)\mapsto \mu\circ (f\stac R g)\circ \Delta,\qquad \textrm{for }f,g\in
{}_R\mathrm{Hom}_R(\cD,A), 
$$
and unit element $\eta\circ \epsilon$.
In parallel to \prref{pr:coalg.cleft}, following \cite[Proposition
  6.4]{BohmVer:Mor&cleft}  
characterises those algebra extensions by an $R$-coring $\cD$ which are
$\cD$-cleft. 

\begin{proposition}
Let $(\cD,\Delta,\epsilon)$ be an $R$-coring and $B\subseteq A$ be an algebra
extension by $\cD$. The $\cD$-extension $B\subseteq A$ is $\cD$-cleft (with
respect to the given $\cD$-comodule structure of $A$) if and only if the
following conditions hold. 

(1) $\cC:=A\otimes_R \cD$ is an $A$-coring, with coproduct $A\otimes_R \Delta$
    and counit $A\otimes_R \epsilon$, the left regular $A$-module structure
    of the first factor and a right $A$-module structure such that the
    $\cD$-coaction in $A$ is a right $A$-module map $A\to \cC$.

(2) $B$ is an $R$-subring of $A$.

(3) There exists a convolution invertible left $R$-module, right
    $\cD$-comodule map $j:\cD\to A$.
\end{proposition}

\thref{thm:Mor_cleft=Gal&norm.b} implies the following theorem. 

\begin{theorem}
An algebra extension $B\subseteq A$ by an $R$-coring $\cD$ is cleft if and only
if it is a $\cD$-Galois extension and the {\em normal basis property} holds,
i.e. $A\cong B\otimes_R \cD$ as left $B$-modules right $\cD$-comodules.
\end{theorem}

The following theorem is a consequence of
\thref{thm:corext.cleft.str.str.thm}. 

\begin{theorem}\thlabel{thm:SSTcoring}
Consider an $R$-coring $\cD$ and a $\cD$-cleft algebra extension $B\subseteq
A$. 
If there exist finite sets, $\{v_i\}\subseteq {}_R\mathrm{Hom}_R(\cD,B)$ and 
$\{d_i\}\subseteq \cD$ satisfying $\sum_i v_i(d_i)=1_B$, then the category of
right $B$-modules is equivalent to the category of right comodules for the
$A$-coring $\cC:=A\otimes_R \cD$, via the induction functor $-\otimes_B
A:\cM_B\to \cM^\cC$.  
\end{theorem}

An algebra extension $B\subseteq A$ by a right $R$-bialgebroid $\hH$ is said
to be {\em $\hH$-cleft} provided that it is a cleft extension by the coring
underlying $\hH$. Via the canonical isomorphism $M\otimes_R H\cong
M\otimes_A A\otimes_R H$, for any right $A$-module $M$, right comodules for
the $A$-coring $A\otimes_R H$ are identified with right modules for the monoid
$A$ in $\cM^\hH$. They are called (right-right) {\em relative Hopf
  modules}. Since the unit element in the algebra $\hH$ is grouplike,
\thref{thm:SSTcoring} implies that, for an $\hH$-cleft algebra extension
$B\subseteq A$, the functor $-\otimes_B \, A$ is an equivalence from the
category of right $B$-modules to the category $\cM^\hH_A$ of right-right
relative Hopf modules. 

\subsection{Cleft extensions by 
\underline{pure}
Hopf algebroids}\selabel{sec:pure.hgd.cleft}

As we have seen in \seref{sec:hgd.Gal}, a Hopf algebroid incorporates two
coring (and bialgebroid) structures. Along the lines in
\seref{sec:coring.cleft}, one can consider cleft extensions by either one of
them. However, it turns out that there is a third, more useful notion of a
cleft extension by a Hopf algebroid. It is more useful in the following
sense. First, it is the notion for which -- analogously to the case of
Hopf algebras -- the total algebra of a Hopf algebroid $\hH$ is an $\hH$-cleft
extension of the base algebra. Second, this definition of a cleft extension by
a Hopf algebroid allows one to extend \thref{thm:cleft=cr.prod} to algebra
extensions by Hopf algebroids. 
{\color{blue}
This definition was proposed in
\cite[\emph{Corrigendum}]{BohmBrz:hgd.cleft}. 

\begin{definition}\delabel{def:hgd.cleft}
Let $\hH$ be any (not necessarily pure) Hopf algebroid, with structure maps
denoted as in \deref{def:hgd}. We say that an $\hH$-extension $B\subseteq A$
is {\em cleft} provided that the following properties are obeyed.

(1) In addition to its $R$-ring structure $\eta_R:R\to A$, $A$ is also an
    $L$-ring, with some unit map $\eta_L:L\to A$.

(2) $B$ is an $L$-subring of $A$.

(3) There exist morphisms $j\in {}_L\mathrm{Hom}^\hH(H,A)$ and ${\widetilde
    j}\in {}_R\mathrm{Hom}_L(H,A)$, such that
$$
\mu_R \circ (j\stac R {\widetilde j})\circ \Delta_R =\eta_L\circ \epsilon_L,
\qquad \textrm{and}\qquad  
\mu_L \circ ({\widetilde j}\stac L j)\circ \Delta_L =\eta_R\circ \epsilon_R,
$$
where $\mu_L$ and $\mu_R$ denote
the multiplications in $A$, as an $L$-ring and as an $R$-ring, respectively. 
The module structures in $H$ are determined by the respective coring structures
and the module structures in $A$ are determined by the respective ring
structures (see the proof of \prref{prop:hgdcleft} below). 
\end{definition}
}

{\color{blue}
We divide our study of cleft extensions by Hopf algebroids into two parts. In
this section we restrict to {\em pure} Hopf algebroids. Under this
restriction, we can apply the methods and results from
\seref{sec:cleft.bicom}. Treatment of the general case, i.e. cleft extensions
by arbitrary Hopf algebroids, requiring some new technics, is left to
next \seref{sec:nonpure.hgd.cleft}. 
}
\smallskip

Consider a Hopf algebroid $\hH$ with constituent left bialgebroid $\hH_L$ over
a base algebra $L$, right bialgebroid $\hH_R$ over a base
algebra $R$, and antipde $S$. 
Take a right $\hH$-comodule algebra $A$ with $\hH_R$-coaction $a\mapsto
a^{[0]}\otimes_R a^{[1]}$ and 
$\hH_L$-coaction $a\mapsto a_{[0]}\otimes_L a_{[1]}$, related via
\equref{eq:hgd_comod}. The comodule algebra $A$ determines an $A$-coring
$\cC_R:=A\otimes_R H$, with $A$-$A$ bimodule structure   
\begin{equation}\eqlabel{eq:AotRH}
a_1(a\stac R h)a_2=a_1a{a_2}^{[0]}\stac R h{a_2}^{[1]},\qquad \textrm{for }
a_1,a_2\in A,\ a\stac R h\in A\stac R H.
\end{equation}
Using the notations in \deref{def:hgd}, the coproduct in $\cC_R$ is $A\otimes_R
\Delta_R$ and the counit is $A\otimes_R \epsilon_R$. 
Obviously, by coassociativity of $\Delta_R$, $\cC_R$ is a $\cC_R$-$\hH_R$
bicomodule, via the left $\cC_R$-coaction provided by the coproduct in $\cC_R$
and right $\hH_R$-coaction $A\otimes_R \Delta_R$. That is to say, the
$R$-coring $(H,\Delta_R,\epsilon_R)$ is a right extension of $\cC_R$. But we
have more: by the Hopf algebroid axiom (ii) in \deref{def:hgd}, $\cC_R$ is
also a $\cC_R$-$\hH_L$ bicomodule, via the left $\cC_R$-coaction provided by
the coproduct in $\cC_R$ and right $\hH_L$-coaction $A\otimes_R
\Delta_L$. Thus also the $L$-coring $(H,\Delta_L,\epsilon_L)$ underlying
$\hH_L$ is a right extension of $\cC_R$. 

{\color{blue}
\begin{lemma}\lelabel{lem:pure}
Consider a Hopf algebroid $\hH$, with structure maps denoted as in
\deref{def:hgd}, and a right $\hH$-comodule algebra $A$. If $\hH$ is a pure
Hopf algebroid then (the coring underlying) $\hH_L$ is a pure extension of
the $A$-coring $\cC_R:=A\ot_R H$.
\end{lemma}

\begin{proof}
By the isomorphism 
\begin{equation}\eqlabel{eq:isom}
W\ot_A \cC_R\cong W\ot_R H, 
\end{equation}
for any right $A$-module $W$, right comodules for $\cC_R$ are identified with
right $A$-modules and right $\hH_R$-comodules $M$, such that the
$\hH_R$-coaction $M\to M\ot_R H$ is a right $A$-module map with respect to the
$A$-action of the codomain arising from the isomorphism \equref{eq:isom}. 

If $\hH$ is a pure Hopf algebroid then in particular the equalizer
\equref{eq:R.M} in $\cM_L$ is $H\ot_L H$-pure, for any $\cC_R$-comodule
$M$. By the isomorphism \equref{eq:isom}, this means  $H\ot_L H$-purity of the
equalizer   
$$
\xymatrix{
M \ar[rr]&&M\stac A \cC_R \ar@<2pt>[rr]\ar@<-2pt>[rr]&& M\stac A \cC_R \stac A 
\cC_R
}
$$
in $\cM_L$, that is, purity of the coring extension in the claim.
\end{proof}
}

Analogously to \prref{prop:Hopf.cleft<Mor}, one proves following
\cite[Proposition 6.6]{BohmVer:Mor&cleft}. It tells us that,
{\color{blue} for a pure Hopf algebroid $\hH$}
and an algebra extension $B\subseteq A$ by $\hH$,
the right $\cC_R$-comodule $A$ extends to a cleft bicomodule for the coring
extension $\hH_L$ of $\cC_R$, introduced above, if and only if $B\subseteq A$
is an $\hH$-cleft extension in the sense of \cite[Definition
  3.5]{BohmBrz:hgd.cleft} (recalled in \deref{def:hgd} below). Recall that a
right comodule  
algebra $A$ of a Hopf algebroid $\hH$ is in particular an $R$-ring, with unit
map $\eta_R:R\to A$, over the base algebra $R$ of the constituent right
bialgebroid in $\hH$. 

\begin{proposition}\prlabel{prop:hgdcleft}
Consider a 
{\color{blue} pure}
Hopf algebroid $\hH$ with constituent left bialgebroid $\hH_L$ over
a base algebra $L$ and right bialgebroid $\hH_R$ over a base
algebra $R$, and an algebra extension $B\subseteq A$ by $\hH$.
The right $\cC_R:=A\ot_R H$-comodule $A$ extends to a cleft bicomodule for the
coring extension $\hH_L$ of $\cC_R$ if and only if $B\subseteq A$ is an
$\hH$-cleft extension.
\end{proposition}
\begin{proof}(Sketch.)
Since $\hH$ is a pure Hopf algebroid, $A$ is a right $\hH$-comodule algebra if
and only if $A$ is a right $\hH_R$-comodule algebra, in which case it is a
right $\cC_R$-comodule as well.
The right $\cC_R$-comodule $A$ extends to
an $L$-$\cC_R$ bicomodule if and only if it is also a left $L$-module such
that the left $L$-action is a right $A$-module map and a right $\hH_R$-comodule
map. The left $L$-action is a right $A$-module map if and only if there exists
an algebra map $\eta_L:L\to A$ in terms of which the left action by $l\in L$
on $A$ is given 
by left multiplication by $\eta_L(l)$. The left $L$-action is a right
$\hH_R$-comodule map if and only if $\eta_L(l)\in B$, for all $l\in L$. Hence
properties (1) and (2) in \deref{def:hgd.cleft} are equivalent to $A$ being an
$L$-$\cC_R$ bicomodule. 

Analogously to \equref{eq:Hopf.Mor}, the Morita
context $\mathbb{M}(A)$, corresponding the the $L$-$\cC_R$ bicomodule $A$ via 
\equref{eq:Mor}, is isomorphic to
\begin{equation}\eqlabel{eq:hgd.Mor}
\big(\ 
{}_L\mathrm{Hom}_L(H,B)\ ,\ 
X\ ,\ 
{}_L\mathrm{Hom}^\hH(H,A)\ ,\ 
{}_{L^{op}}\mathrm{Hom}^\hH(H^{tw},A)\ ,\ 
\blacklozenge'\ ,\  
\lozenge'\big).
\end{equation}
The (co)module structures in \equref{eq:hgd.Mor} need to be explained. 
Let us use the notations introduced in, and after \deref{def:hgd}, and in the 
paragraph after \deref{def:hgd_comod}.
An element $v\in {}_L\mathrm{Hom}_L(H,B)$ is a
bimodule map with respect to the actions
\begin{equation}\eqlabel{eq:bim.L-L}
v(s_L(l) t_L(l') h)=\eta_L(l)v(h)\eta_L(l'),\qquad \textrm{for } l,l'\in L,\
h\in H.
\end{equation}
As a $k$-module, 
\begin{equation}\eqlabel{eq:X}
X
=
\{\ {\widetilde u}\in
{}_R \mathrm{Hom}_R(H,A)\ |\ {\widetilde u}(h^{(2)})^{[0]}
\stac R h^{(1)}{\widetilde u}(h^{(2)})^{[1]}={\widetilde u}(h^{(1)})\stac R
h^{(2)},\quad \forall h\in H\ \},
\end{equation}
where ${\widetilde u}\in {}_R \mathrm{Hom}_R(H,A)$ is a bimodule map with
respect to the actions 
\begin{equation}\eqlabel{eq:bim.R-R}
{\widetilde u}(h s_R(r) t_R(r'))=\eta_R(r'){\widetilde u}(h)\eta_R(r),\qquad
\textrm{for } r,r'\in R,\ h\in H.
\end{equation}
In ${}_L\mathrm{Hom}^\hH(H,A)$, $H$ is a left $L$-module via $s_L$, 
and $A$ is a left $L$-module via $\eta_L$. Note that the
coproducts in $\hH$ are left $L$-module maps and (since $B$ is an $L$-ring by
definition) so are the $\hH$-coactions in $A$. Thus both $H$ and $A$ are
$L$-$\hH$ bicomodules. Since an element $p\in {}_L\mathrm{Hom}^\hH(H,A)$ is
a comodule map with respect to the regular $\hH$-coactions, it is $R$-$R$
bilinear in the sense that 
$p(s_R(r)hs_R(r'))= \eta_R(r)p(h)\eta_R(r')$, for $h\in H$ and $r,r'\in R$. In
particular, $p$ is an $L$-$R$ bimodule map with respect to the actions
\begin{equation}\eqlabel{eq:bim.L-R}
 p(s_L(l)hs_R(r))=\eta_L(l)p(h)\eta_R(r),\quad \textrm{for } h\in H,\  l\in
 L,\ r\in R.
\end{equation}
In ${}_{L^{op}}\mathrm{Hom}^\hH(H^{tw},A)$, $H^{tw}$ is the same $k$-module
$H$. It is considered to be a left $L^{op}$-module via $t_L$,
and $A$ is a left $L^{op}$-module via right multiplication
by $\eta_L$. $H^{tw}$ is a right $R$-module via $t_R$, a right $L$-module via
$t_R \circ \varepsilon_R \circ t_L$, and a right $\hH$-comodule via the
twisted $\hH_R$-coaction $h\mapsto h_{(2)}\otimes_R S(h_{(1)})$ and 
$\hH_L$-coaction $h\mapsto h^{(2)}\otimes_L S(h^{(1)})$. 
Note that both the twisted coactions in $H^{tw}$ and the $\hH$-coactions in $A$
are left $L^{op}$-module maps, thus both $H^{tw}$ and $A$ are $L^{op}$-$\hH$
bicomodules. An element $q\in {}_{L^{op}}\mathrm{Hom}^\hH(H^{tw},A)$ is
$R$-$R$ bilinear in the sense that 
$q(t_R(r) h t_R(r'))= \eta_R(r')q(h)\eta_R(r)$, for $h\in H$ and $r,r'\in
R$. In particular, $q$ is an $R$-$L$ bimodule map with respect to the actions
\begin{equation}\eqlabel{eq:bim.R-L}
q(t_L(l)h t_R(r))=\eta_R(r)q(h)\eta_L(l),\quad \textrm{for }h\in H,\  l\in L,\ 
r\in R. 
\end{equation}
With the bimodule structures \equref{eq:bim.L-L}, \equref{eq:bim.R-R},
\equref{eq:bim.L-R} and \equref{eq:bim.R-L} in mind,
in the Morita context \equref{eq:hgd.Mor} all algebra and bimodule structures
and also the connecting maps 
$\blacklozenge'$ and $\lozenge'$ can be written as a generalized convolution
product in \cite[Section 3]{BohmBrz:hgd.cleft}: 
$$
(f,g)\mapsto  \mu_Q\circ (f\stac Q g)\circ \Delta_Q,\qquad \textrm{for
}P,Q,S\in \{L,R\},\ f\in {}_P\mathrm{Hom}_Q(H,A),\ g\in
{}_Q\mathrm{Hom}_S(H,A),
$$
where $\mu_Q$ denotes multiplication in (the $Q$-ring) $A$.
Thus the equivalence of the cleft property of the $L$-$\cC_R$ bicomodule $A$
and property (3) in \deref{def:hgd.cleft} follows by \cite[Lemmata 3.7 and
  3.8]{BohmBrz:hgd.cleft}, stating that the (generalized) convolution inverse
of $j\in {}_L\mathrm{Hom}^\hH(H,A)$, if it exists, belongs to
${}_{L^{op}}\mathrm{Hom}^\hH(H^{tw},A)$.  
\end{proof}

In the same way as \deref{def:hgd.cleft} of a cleft extension by a Hopf
algebroid $\hH$  
combines the two bialgebroids in $\hH$, so does the following analogue of
\thref{thm:cleft=Gal&norm.b}, which is a consequence of
\thref{thm:Mor_cleft=Gal&norm.b}. It was proven for arbitrary Hopf algebroids
in \cite[Theorem 3.12]{BohmBrz:hgd.cleft} (see \thref{thm:nonpure.cleft}).

\begin{theorem} \thlabel{thm:pure_cleft}
Consider a 
{\color{blue} pure}
Hopf algebroid $\hH$ with constituent left bialgebroid $\hH_L$ over
a base algebra $L$ and right bialgebroid $\hH_R$ over a base algebra $R$. An
algebra extension $B\subseteq A$ by $\hH$ is cleft if and only if the
following properties hold. 

(1) The extension $B\subseteq A$ is $\hH_R$-Galois.

(2) The {\em normal basis property} holds, i.e. $A\cong B\otimes_L H$, as left
    $B$-modules and right $\hH$-comodules.
\end{theorem}
Since the unit element of the algebra underlying a Hopf algebroid $\hH$ is
grouplike, \thref{thm:corext.cleft.str.str.thm} implies the following Strong
Structure Theorem. 
Recall that, for a right $\hH$-comodule algebra $A$, a right-right
$(\hH,A)$-relative Hopf module is defined as a right module of the monoid $A$
in $\cM^\hH$. 
{\color{blue} If $\hH$ is a pure Hopf algebroid then} an $(\hH,A)$-relative
Hopf module is is canonically identified with a right comodule for the
$A$-coring \equref{eq:AotRH}. 

\begin{theorem}\thlabel{thm:hgd.cleft.str.str.thm}
For a cleft algebra extension $B\subseteq A$ by a 
{\color{blue} pure}
Hopf algebroid $\hH$, the category of
right $B$-modules is equivalent to the category $\cM^\hH_A$ of right-right
relative Hopf modules, via the induction functor $-\otimes_B
A:\cM_B\to \cM^\hH_A$.  
\end{theorem}
A generalization of \thref{thm:hgd.cleft.str.str.thm}
for any Hopf algebroid is given in \thref{thm:nonpure.str.str.thm}.

\subsection{Cleft extensions by arbitrary Hopf
  algebroids}\selabel{sec:nonpure.hgd.cleft} 

{\color{blue}
In \seref{sec:pure.hgd.cleft} we applied the theory of cleft bicomodules of
pure coring extensions to study cleft extensions by {\em pure} Hopf algebroids.
In \deref{def:hgd.cleft},
the cleft property of algebra extensions by
arbitrary Hopf algebroids was defined. 
The aim of the current section is to show that also the definition of any Hopf
algebroid cleft extension in \deref{def:hgd.cleft}
is equivalent to the existence of mutually inverse
elements in the Morita context \equref{eq:hgd.Mor}, which is perfectly
meaningful for any Hopf algebroid. 
Although derivation of the the Morita context \equref{eq:hgd.Mor} in the
general case below is very similar to the methods in the previous sections, it
is not known to correspond to any coring extension.
In the final part of the section we describe cleft extensions of arbitrary
Hopf algebroids as crossed products with invertible cocycles, extending in
this way \thref{thm:cleft=cr.prod}. 
\bigskip

For an arbitrary Hopf algebroid $\hH$ and a right $\hH$-comodule algebra $A$,
the coring extension in \leref{lem:pure} is not known to be pure. Hence a
cotensor-functor $-\Box_{\cC_R}\, \cC_R:\cM^{\cC_R} \to \cM^{\hH_L}$ is not
available to derive the Morita context used in the previous section. 
However, if we assume that there is an algebra map $\eta_L$ from the base
algebra $L$ of the constituent left bialgebroid $\hH_L$  to
$B:=A^{co\hH_R}\subseteq A$, then we are equipped with two functors from the
category $\cM^\hH_A$ of right $A$-modules in $\cM^\hH$ to $\cM^\hH$: The
forgetful functor $U$ and $V:=(-)^{co\hH_R}\ot_L H$. (Recall from
the proof of \coref{cor:hgd.coinv} that 
$(-)^{co\hH_R}\cong \mathrm{Hom}^{\hH_R}_A(A,-)=\mathrm{Hom}^\hH_A(A,-)$
and, since for $M\in \cM^\hH_A$, $M^{co\hH_R}$ is a 
right $B$-module and $B$ is an $L$-ring, $M^{co\hH_R}$ is a right $L$-module.)
Clearly, if $\hH$ is a pure Hopf algebroid then $U$ and $V$ differ by trivial
isomorphisms from the respective functors $-\Box_{\cC_R}\, \cC_R$ and
$\mathrm{Hom}^{\cC_R}(A,-)\ot_L H:\cM^{\cC_R}\to \cM^{\hH_L}$.
The following proposition is a detailed version of an observation in
\cite[\emph{Corrigendum}]{BohmVer:Mor&cleft}. 

\begin{proposition}\prlabel{pr:iso_mor}
Let $\hH$ be a Hopf algebroid with structure maps denoted as in
\deref{def:hgd} and $B\subseteq A$ be an $\hH$-extension. Assume that there
exists an algebra map $\eta_L:L\to B=A^{co\hH_R}$ and consider the above
functors $U$ and $V$ from $\cM^\hH_A$ to $\cM^\hH$. The corresponding Morita
context \equref{eq:cat_morita} is isomorphic to the Morita context
\equref{eq:hgd.Mor}. 
\end{proposition}

\begin{proof}
The four maps establishing the stated isomorphism of Morita contexts are very
similar to those in the proof of \cite[Proposition 3.1]{BohmVer:Mor&cleft}:
$$
\begin{array}{lll}
\alpha_1&:{}_L\mathrm{Hom}_L(H,B) \to \mathrm{Nat}(V,V),\qquad 
&v \mapsto [\Phi^v_M: M^{co\hH_R}\ot_L H\to M^{co\hH_R}\ot_L H,\\
&&\hspace{2.85cm} n\ot_L h \mapsto nv(h_{(1)})\ot_L h_{(2)}];\\
\alpha_1^{-1}&:\mathrm{Nat}(V,V) \to {}_L\mathrm{Hom}_L(H,B),\qquad
&\Phi \mapsto (B\ot_L \varepsilon_L)\big(\Phi_A(1_B\ot_L -)\big);\\
\alpha_2&:X\to \mathrm{Nat}(U,U),\qquad \qquad \qquad 
&u \mapsto [\Xi^u_M: M\to M,\ m\mapsto m^{[0]}u(m^{[1]})];\\
\alpha_2^{-1}&:\mathrm{Nat}(U,U) \to X\qquad \qquad \qquad 
&\Xi\mapsto (A\ot_R \varepsilon_R)\big(\Xi_{A\ot_R H}(1_A \ot_R -)\big);\\
\alpha_3&:{}_L\mathrm{Hom}^\hH(H,A)\to \mathrm{Nat}(V,U),\qquad 
&p\mapsto [\Theta^p_M:M^{co\hH_R}\ot_L H\to M,\ 
n\ot_L h \mapsto n p(h)];\\
\alpha_3^{-1}&:\mathrm{Nat}(V,U) \to {}_L\mathrm{Hom}^\hH(H,A),\qquad
&\Theta \mapsto \Theta_A(1_B \ot_L -);\\
\alpha_4&:{}_{L^{op}}\mathrm{Hom}^\hH(H^{tw},A)\to \mathrm{Nat}(U,V),\qquad 
&q\mapsto [\Omega^q_M:M \to M^{co\hH_R}\ot_L H,\\
&& \hspace{1.9cm}m \mapsto {m_{[0]}}^{[0]}q({m_{[0]}}^{[1]})\ot_L m_{[1]}];\\
\alpha_4^{-1}&:\mathrm{Nat}(U,V) \to {}_{L^{op}}\mathrm{Hom}^\hH(H^{tw},A),\qquad
&\Omega \mapsto (A\ot_L \varepsilon_L)\big(\Omega_{A\ot_R H}(1_A\ot_R -)\big).\\
\end{array}
$$
By the $L$-$L$ bilinearity of $v$, $\Phi^v_M$ is a well defined map
$M^{co\hH_R}\ot_L H \to M \ot_L H$. Since $M^{co\hH_R}$ is a right $B$-module
and the range of $v$ is in $B$, the range of $\Phi^v_M$ is in
$M^{co\hH_R}\ot_L H$ as needed. By the right $\hH$-colinearity of the
coproduct $\Delta_L$ in $\hH_L$, $\Phi^v_M$ is
right $\hH$-colinear. Naturality of $\Phi^v$ is obvious. In order to see that 
$\alpha_1$ is a ring map, compute for $n\ot_L h \in M^{co\hH_R}\ot_L H$ and
$v,v'\in {}_L \mathrm{Hom}_L(H,B)$,
\begin{eqnarray*}
\Phi^{\eta_L\circ \varepsilon_L}_M(n\stac L h)
&=& n \eta_L(\varepsilon_L(h_{(1)})) \stac L h_{(2)} 
= n\stac L s_L(h_{(1)})) h_{(2)} 
= n \stac L h;\\
\big(\Phi^v_M\circ \Phi^{v'}_M)(n\stac L h)
&=& n v'(h_{(1)}) v(h_{(2)}) \stac L h_{(3)}
= n (v'  v)(h_{(1)}) \stac L h_{(2)}
= \Phi^{v'  v}_M(n\stac L h).
\end{eqnarray*}
In the converse direction, for $\Phi\in \mathrm{Nat}(V,V)$, $\Phi_A$ is right
$\hH$-colinear hence in particular right $L$-linear. Since also the
$\hH_L$-counit $\varepsilon_L$ is right $L$-linear, so is
$\alpha_1^{-1}(\Phi)$. In order to see that $\alpha_1^{-1}(\Phi)$ is also left
$L$-linear, note that for any $M\in\cM^\hH_A$ and $n\in M^{co\hH_R}$, the map
\begin{equation}\eqlabel{eq:n.map}
A\to M,\qquad a\mapsto na
\end{equation}
is a morphism in $\cM^\hH_A$. Thus by naturality of 
$\Phi$, 
\begin{equation}\eqlabel{eq:n_4}
\Phi_M(nb\stac L h)=n \Phi_A (b\stac L h),\qquad \textrm{for }n\in
M^{co\hH_R},\ b\ot_L h\in  B\ot_L H.
\end{equation}
Applying it to $M=A$, we conclude that $\Phi_A$ is left $B$-linear thus in
particular left $L$-linear. 
This proves the left $L$-linearity of $\alpha_1^{-1}(\Phi)$.
By the right $L$-linearity of $v\in
{}_L\mathrm{Hom}_L(H,B)$ and the counitality of $\Delta_L$, for all $h\in H$, 
$$
(\alpha_1^{-1}\circ \alpha_1)(v)(h)
=v(h_{(1)})\varepsilon_L(h_{(2)})=v(h).
$$
In order to see that $\alpha_1^{-1}$ is also the right inverse of $\alpha_1$,
use the right $\hH_L$-colinearity of $\Phi_A$, i.e. the identity
$(B\ot_L \Delta_L)\circ \Phi_A =(\Phi_A \ot_L H)\circ (B\ot_L \Delta_L)$,
(in the second equality) and \equref{eq:n_4} (in the
last equality), for $n\ot_L h\in M^{co\hH_R}\ot_L H$:
$$
(\alpha_1\circ \alpha_1^{-1})(\Phi)(n\stac L h)
= n (B\stac L \varepsilon_L)\big(\Phi_A(1_B\stac L h_{(1)})\big)\stac L h_{(2)}
= n \Phi_A(1_B\stac L h) =\Phi_M(n\stac L h).
$$

By the left $R$-linearity of $u\in X$, $\Xi^u_M$ is a well defined map $M\to
M$. In order to check the $\hH_R$-colinearity of $\Xi^u_M$, use the
$\hH_R$-colinearity of the $A$-action on $M$ and the defining condition
\equref{eq:X} of $X$:
$$
\Xi^u_M(m)^{[0]} \stac R \Xi^u_M(m)^{[1]}
= m^{[0]} u(m^{[2]})^{[0]} \stac R m^{[1]}u(m^{[2]})^{[1]}=
m^{[0]} u(m^{[1]})\stac R m^{[2]}
= \Xi^u_M(m^{[0]})\stac R m^{[1]},
$$
for all $m\in M$. We conclude by \prref{prop:ff} that $\Xi^u_M$ is right
$\hH$-colinear. Naturality of $\Xi^u$ is obvious. In order to check that
$\alpha_2$ is a ring map, compute for $m\in M$, $u,u'\in X$,
\begin{eqnarray*}
\Xi^{\eta_R \circ \varepsilon_R}_M(m)
&=& m^{[0]} \eta_R(\varepsilon_R(m^{[1]}))
=m;\\
(\Xi^u_M \circ \Xi^{u'}_M)(m)
&=& m^{[0]} u'(m^{[2]})^{[0]} u\big(m^{[1]} u'(m^{[2]})^{[1]}\big)\\
&=& m^{[0]} u'(m^{[1]})u(m^{[2]})
= m^{[0]} (u' u)(m^{[1]})
= \Xi^{u' u}_M(m),
\end{eqnarray*}
where in the first equality of the last line we used \equref{eq:X}.

In the converse direction, for $\Xi\in  \mathrm{Nat}(U,U)$, $\Xi_{A\ot_R H}$
is right $\hH$-colinear hence in particular right $R$-linear.
Since also the $\hH_R$-counit $\varepsilon_R$ is right $R$-linear, so is
$\alpha_2^{-1}(\Xi)$. 
For any right $A$-module $W$, $W\ot_R H$ is an object in $\cM^\hH_A$, via the
$\hH$-comodule structure of the second factor (both coproducts $\Delta_L$ and
$\Delta_R$ are left $R$-linear) and the right $A$-action $(w\ot_R h)a=
wa^{[0]} \ot_R ha^{[1]}$. Moreover, for any $w\in W$, the map 
\begin{equation}\eqlabel{eq:p.map}
A \stac R H \to W\stac R H,\qquad a\stac R h \mapsto wa\stac R h
\end{equation}
is a morphism in $\cM^\hH_A$. Thus by the naturality of $\Xi$,
\begin{equation}\eqlabel{eq:n_5}
\Xi_{W\ot_R H}(wa\stac R h) = w \Xi_{A\ot_R H}(a\stac R h),\quad 
\textrm{for }w\in W,\ a\stac R h\in A\stac R H.
\end{equation}
Applying it to $W=A$, we conclude that $\Xi_{A\ot_R H}$ is left $A$-linear
hence in particular left $R$-linear. This proves the left $R$-linearity of
$\alpha_2^{-1}(\Xi)$.  
In order to see that $\alpha_2^{-1}(\Xi)$ belongs to $X$, observe that
for any $M\in \cM^\hH_A$, the $\hH_R$-coaction $M\to M\ot_R H$, $m\mapsto
m^{[0]}\ot_R m^{[1]}$ is a morphism in $\cM^\hH_A$. Thus by the naturality of
$\Xi$,
\begin{equation}\eqlabel{eq:n_7}
\Xi_M(m)^{[0]}\stac R \Xi_M(m)^{[1]}= \Xi_{M\ot_R H}(m^{[0]}\stac R
m^{[1]}),\qquad \textrm{for }m\in M.
\end{equation}
Applying \equref{eq:n_5} for $W=A\ot_R H$ (in the second equality),
\equref{eq:n_7} for $M=A\ot_R H$ (in the third equality) and the right
$\hH_R$-colinearity of $\Xi_{A\ot_R H}$ (in the last equality), one computes
for $h\in H$, 
\begin{eqnarray*}
\alpha_2^{-1}(\Xi)(h^{(2)})^{[0]} &\stac R &h^{(1)}
\alpha_2^{-1}(\Xi)(h^{(2)})^{[1]} 
= (1_A \stac R h^{(1)})
(A\stac R \varepsilon_R)\big(\Xi_{A\ot_R H}(1_A \stac R h^{(2)})\big)\\
&=& (A\stac R H \stac R \varepsilon_R)\big(
\Xi_{A\ot_R H \ot_R H}(1_A \stac R h^{(1)}\stac R h^{(2)})\big)\\
&=& \big((A\stac R H\stac R \varepsilon_R)\circ (A\stac R \Delta_R)\big)
\big(\Xi_{A\ot_R  H}(1_A \stac R h)\big)
= \Xi_{A\ot_R  H}(1_A \stac R h)\\
&=& \big((A\stac R \varepsilon_R\stac R H)\circ (A\stac R \Delta_R)\big)
\big(\Xi_{A\ot_R  H}(1_A \stac R h)\big)\\ 
&=&\alpha_2^{-1}(\Xi)(1_A \stac R h^{(1)}) \stac R h^{(2)}.
\end{eqnarray*}
Making use of \equref{eq:X}, the right $R$-linearity of $u\in X$ and the
counitality of $\Delta_R$, we see that for all $h\in H$,
$$
(\alpha_2^{-1} \circ\alpha_2)(u)(h) 
= u(h^{(2)})^{[0]} \eta_R\big(\varepsilon_R(h^{(1)}u(h^{(2)})^{[1]})\big)
= u(h^{(1)})\eta_R(\varepsilon_R(h^{(2)})) = u(h).
$$
Using \equref{eq:n_5} (in the second equality) and \equref{eq:n_7} (in the
third equality) one checks that for $\Xi\in \mathrm{Nat}(U,U)$ and $m\in M$,
also 
\begin{eqnarray*}
(\alpha_2 \circ \alpha_2^{-1})(\Xi)(m)
&=&(M\stac R \varepsilon_R)\big(m^{[0]}\Xi_{A\ot_R H}(1_A \stac R m^{[1]})\big)
=(M\stac R \varepsilon_R)\big(\Xi_{M\ot_R H}(m^{[0]}\stac R m^{[1]})\big)\\
&=& \Xi_M(m)^{[0]}\eta_R(\varepsilon_R(\Xi_M(m)^{[1]}))
=\Xi_M(m).
\end{eqnarray*}

By the left $L$-linearity of $p\in {}_L\mathrm{Hom}^\hH(H,A)$, $\Theta^p_M$ is
a well defined map. By the right $\hH$-colinearity of $p$, $\Theta^p_M$ is
$\hH$-colinear. Naturality of $\Theta$ is obvious. Compatibility of $\alpha_3$
with the bimodule structures is checked by the following simple computations,
for $p\in {}_L\mathrm{Hom}^\hH(H,A)$, $v\in {}_L\mathrm{Hom}_L(H,B)$, $u\in
X$ and $n\ot_L h\in M^{co\hH_R}\ot_L H$.
\begin{eqnarray*}
(\Theta^p_M\circ \Phi^v_M)(n\stac L h)
&=&n v(h_{(1)}) p(h_{(2)}) 
= n (vp)(h)
= \Theta^{vp}_M(n\stac L h);\\
(\Xi^u_M\circ \Theta^p_M)(n\stac L h)
&=& n p(h)^{[0]} u\big(p(h)^{[1]}\big)
= n p(h^{(1)}) u(h^{(2)})
= n (pu)(h)
= \Theta^{pu}_M(n\stac L h).
\end{eqnarray*}

In the converse direction, using that \equref{eq:n.map} is a morphism in
$\cM^\hH_A$, the naturality of $\Theta\in \mathrm{Nat}(V,U)$ implies that
\begin{equation}\eqlabel{eq:n_8}
\Theta_M(nb\stac L h)= n \Theta_A(b\stac L h),\qquad \textrm{ for }n\in
M^{co\hH_R},\, b\stac L h\in B\stac L H.
\end{equation}
In particular, $\Theta_A$ is left $B$-linear hence in particular left
$L$-linear. This proves that $\alpha_3^{-1}(\Theta)$ is left $L$-linear. 
By the right $\hH$-colinearity of $\Theta_A$, $\alpha_3^{-1}(\Theta)$ is right
$\hH$-colinear. It is an immediate consequence of the construction of
$\alpha_3$ and $\alpha_3^{-1}$ that $(\alpha_3^{-1}\circ \alpha_3)(p)=p$, for
$p\in {}_L\mathrm{Hom}^\hH(H,A)$. 
It follows easily by \equref{eq:n_8} that also $(\alpha_3 \circ
\alpha_3^{-1})(\Theta) = \Theta$, for $\Theta\in \mathrm{Nat}(V,U)$.

By the left $L^{op}$-, and the left $R$-linearity of $q\in
{}_{L^{op}}\mathrm{Hom}^\hH(H^{tw},A)$, $\Omega^q_M$ is a well defined map $M
\to M\ot_L H$. In order to see that the range of $\Omega^q_M$ is in
$M^{co\hH_R}\ot_L H$, as needed, check that $m^{[0]}q(m^{[1]})\in M^{co\hH_R}$,
for any $M\in \cM^\hH_A$ and $m\in M$:
\begin{eqnarray}
\big(m^{[0]}q(m^{[1]})\big)^{[0]}&\stac R &\big(m^{[0]}q(m^{[1]})\big)^{[1]}
= m^{[0]} q(m^{[2]})^{[0]} \stac R m^{[1]} q(m^{[2]})^{[1]}\nonumber\\
&=& m^{[0]} q({m^{[2]}}_{(2)})\stac R m^{[1]} S({m^{[2]}}_{(1)})
= m^{[0]} q({m^{[1]}}_{(2)})\stac R {{m^{[1]}}_{(1)}}^{(1)}
S({{m^{[1]}}_{(1)}}^{(2)}) \nonumber\\
&=& m^{[0]} q({m^{[1]}}_{(2)})\stac R
s_L(\varepsilon_L({m^{[1]}}_{(1)}))\nonumber \\
&=& m^{[0]} q({m^{[1]}}_{(2)})
\eta_R(\varepsilon_R(s_L(\varepsilon_L({m^{[1]}}_{(1)})))) \stac R 1_H
= m^{[0]}q(m^{[1]})\stac R 1_H.\eqlabel{eq:q_coinv}
\end{eqnarray} 
In the first equality we used the right $\hH_R$-colinearity of the $A$-action
on $M$.
In the second equality we used the right $\hH_R$-colinearity of $q$.
In the third equality we used the right $\hH_L$-colinearity of the
$\hH_R$-coproduct $\Delta_R$. 
The fourth equality follows by the second one of the antipode axioms (iv) in
\deref{def:hgd}. 
In the penultimate equality we applied the Hopf algebroid axiom $t_R\circ
\varepsilon_R \circ s_L =s_L$, cf. \deref{def:hgd} (i).
The same Hopf algebroid axiom $t_R\circ \varepsilon_R \circ s_L =s_L$ is
applied again in the last equality, together with 
the right $R$-linearity of $q$ and the counitality of $\Delta_L$.
By the right $\hH$-colinearity of the $\hH_L$-coaction on $M$, $\Omega^q_M$ is
right $\hH$-colinear. Naturality of $\Omega^q$ is obvious.
The compatibility of $\alpha_4$ with the bimodule structures is checked by the
following computations, for $v\in {}_L\mathrm{Hom}_L(H,B)$, $u\in X$, $q\in
{}_{L^{op}}\mathrm{Hom}^\hH(H^{tw},A)$ and $m\in M$:
\begin{eqnarray*}
(\Phi^v_M\circ \Omega^q_M)(m)
&=& {m_{[0]}}^{[0]}q({m_{[0]}}^{[1]}) v(m_{[1]}) \stac L m_{[2]}
= m^{[0]}q({m^{[1]}}_{(1)}) v({m^{[1]}}_{(2)})\stac L {m^{[1]}}_{(3)}\\
&=& m^{[0]} (qv)({m^{[1]}}_{(1)})\stac L {m^{[1]}}_{(2)}
= {m_{[0]}}^{[0]}(qv)({m_{[0]}}^{[1]})\stac L m_{[1]}
= \Omega^{qv}_M(m);\\
(\Omega^q_M\circ \Xi^u_M)(m)
&=& {\big(m^{[0]} u(m^{[1]})\big)_{[0]}}^{[0]}
q\left( {\big(m^{[0]} u(m^{[1]})\big)_{[0]}}^{[1]}\right) \stac L
\big(m^{[0]} u(m^{[1]})\big)_{[1]}\\
&=& {\big(m^{[0]} u(m^{[1]})\big)}^{[0]}
q\left( {\big(m^{[0]} u(m^{[1]})\big)^{[1]}}_{(1)}\right) \stac L
{\big(m^{[0]} u(m^{[1]})\big)^{[1]}}_{(2)}\\
&=& m^{[0]} u(m^{[2]})^{[0]}
q\big(({m^{[1]}}{u(m^{[2]})^{[1]}})_{(1)}\big) \stac L 
({m^{[1]}}{u(m^{[2]})^{[1]}})_{(2)}\\
&=&  m^{[0]}u(m^{[1]}) q({m^{[2]}}_{(1)}) \stac L {m^{[2]}}_{(2)}
= m^{[0]} u({{m^{[1]}}_{(1)}}^{(1)}) q({{m^{[1]}}_{(1)}}^{(2)}) \stac L
{m^{[1]}}_{(2)} \\
&=&  m^{[0]} (uq)({m^{[1]}}_{(1)})\stac L {m^{[1]}}_{(2)} 
= {m_{[0]}}^{[0]} (uq)({m_{[0]}}^{[1]}) \stac L m_{[1]}
=\Omega^{uq}_M(m).
\end{eqnarray*}
In the second and penultimate equalities of both computations we used the right
$\hH_L$-colinearity of the right $\hH_R$-coaction on $M$. 
In the third equality of the second computation we used the right
$\hH_R$-colinearity of the $A$-action on $M$ and coassociativity of the 
$\hH_R$-coaction. 
In the fourth equality \equref{eq:X} has been used and in the fifth equality we
made use of the right $\hH_L$-colinearity of $\Delta_R$. 

In the definition of $\alpha_4^{-1}$ the isomorphism
\begin{equation}\eqlabel{eq:co_iso}
W \to (W\stac R H)^{co\hH_R},\ w\mapsto w\stac R 1_H;\qquad 
(W\stac R H)^{co\hH_R} \to W,\ \sum_i w_i \stac R h_i \mapsto \sum_i w_i
\eta_R\big(\varepsilon_R(h_i)\big) 
\end{equation}
is implicitly used, for $W\in \cM_A$. This isomorphism induces a right
$B$-action on $(W\stac R H)^{co\hH_R}$, as $(\sum_i w_i \stac R h_i)b=\sum_i
w_ib \stac R h_i$ (meaningful in light of \reref{rem:B_R}). In particular,
since $B$ is an $L$-ring,
\equref{eq:co_iso} is a right $L$-module isomorphism. For $\Omega\in
\mathrm{Nat}(U,V)$, the map 
$\Omega_{A\ot_R H}$ is right $\hH$-colinear hence in particular right
$L$-linear. Since also $\varepsilon_L$ is right $L$-linear, so is
$\alpha_4^{-1} (\Omega)$. For any $r\in R$, the map $A\ot_R H \to A\ot_R H$,
$a\ot_R h\mapsto a\ot_R t_R(r)h$ is a morphism in $\cM^\hH_A$. Hence it
follows by the naturality of $\Omega$ that $\alpha_4^{-1} (\Omega)$ is right
$R$-linear. 
Since the right $\hH_R$-coaction on any $M\in \cM^\hH_A$ is a morphism in
$\cM^\hH_A$, naturality of $\Omega$ implies
\begin{equation}\eqlabel{eq:n_13}
\Omega_{M\ot_R H}(m^{[0]}\stac R m^{[1]})=\Omega_M(m)\in M\stac L H, \qquad
\textrm{for }m\in M,
\end{equation}
where the isomorphism \equref{eq:co_iso} is implicitly used.
Applying $M\ot_L \varepsilon_L$ to both sides and using the isomorphism
\equref{eq:co_iso}, it follows that 
\begin{equation}\eqlabel{eq:n_10}
\big( (M\stac R H)^{co\hH_R}\stac L \varepsilon_L\big)
\big(\Omega_{M\ot_R H}(m^{[0]}\stac R m^{[1]})\big)=
(M\stac L \varepsilon_L)\big(\Omega_M(m)\big)\stac R 1_H, 
\quad \textrm{for }m\in M,
\end{equation}
as elements of $(M\ot_R H)^{co\hH_H}\subseteq M\stac R H$.
Moreover, using that \equref{eq:p.map} is a morphism in $\cM^\hH_A$ and the
naturality of $\Omega$, it follows that for any right $A$-module $W$, 
\begin{equation}\eqlabel{eq:n_11}
\Omega_{W\ot_R H}(wa \stac R h)=w \Omega_{A\ot_R H}(a\stac R h),\qquad
\textrm{for }w\in W,\ a\stac R h\in A\stac R H.
\end{equation}
Combining \equref{eq:n_10} for $M=A\ot_R H$, and \equref{eq:n_11} for
$W=A\ot_R H$, we obtain that \begin{equation}\eqlabel{eq:n_12}
\alpha_4^{-1}(\Omega)(h)\stac R 1_H = \alpha_4^{-1}(\Omega)(h^{(2)})^{[0]}
\stac R h^{(1)}  \alpha_4^{-1}(\Omega)(h^{(2)})^{[1]}.
\end{equation}
By the right $R$-linearity of $\alpha_4^{-1}(\Omega)$ and the Hopf algebroid
axiom $t_R\circ \varepsilon_R \circ s_L =s_L$, it follows that for $l\in L$
and $h\in H$, 
$\alpha_4^{-1}(\Omega)(s_L(l)h)\ot_R 1_H = \alpha_4^{-1}(\Omega)(h)\ot_R
s_L(l)$. Thus the expression $\alpha_4^{-1}(\Omega)(h_{(2)})\ot_L S(h_{(1)})$
is meaningful, and with identity \equref{eq:n_12} at hand, it satisfies
\begin{eqnarray*}
\alpha_4^{-1}(\Omega)(h_{(2)})\stac R S(h_{(1)})
&=& \alpha_4^{-1}(\Omega)({h_{(2)}}^{(2)})^{[0]} \stac R 
S(h_{(1)}) {h_{(2)}}^{(1)}  \alpha_4^{-1}(\Omega)({h_{(2)}}^{(2)})^{[1]}\\
&=& \alpha_4^{-1}(\Omega)({h}^{(2)})^{[0]}\stac R S({h^{(1)}}_{(1)})
{h^{(1)}}_{(2)}  \alpha_4^{-1}(\Omega)(h^{(2)})^{[1]}\\
&=& \alpha_4^{-1}(\Omega)({h}^{(2)})^{[0]}\stac R s_R(\varepsilon_R(h^{(1)}))
\alpha_4^{-1}(\Omega)(h^{(2)})^{[1]}\\
&=& \big( \eta_R(\varepsilon_R(h^{(1)}))\alpha_4^{-1}
(\Omega)({h}^{(2)})\big)^{[0]} \stac R  
\big( \eta_R(\varepsilon_R(h^{(1)})) \alpha_4^{-1}
(\Omega)({h}^{(2)})\big)^{[1]}\\ 
&=& \alpha_4^{-1}(\Omega)(h)^{[0]}\stac R \alpha_4^{-1}(\Omega)(h)^{[1]}.
\end{eqnarray*}
In the first equality we used \equref{eq:n_12} and in the second one we used
the right $\hH_R$-colinearity of $\Delta_L$. The third equality follows by the
first one of the antipode axioms in \deref{def:hgd}(iv). In the penultimate
equality we used that by the right $R$-linearity and unitality the
$\hH_R$-coaction on $A$, it follows that $\eta_R(r)^{[0]}\ot_R \eta_R(r)^{[1]}
= 1_A \ot_R s_R(r)$, for any $r\in R$, and that the multiplication in $A$ is
a right $\hH_R$-comodule map. The last equality follows by noting
that by \equref{eq:n_11} the map $\Omega_{A\ot_R H}$ is left $A$-linear, hence
$\alpha_4^{-1}(\Omega)$ is left $R$-linear in the sense that, for $r\in 
R$ and $h\in H$,
\begin{eqnarray*}
\eta_R(r)\alpha_4^{-1}(\Omega)(h)
&=& (A\stac L \varepsilon_L)\big(\Omega_{A\ot_R H}(\eta_R(r)\stac R h)\big) \\
&=& (A\stac L \varepsilon_L)\big(\Omega_{A\ot_R H}(1_A \stac R h t_R(r))\big)
= \alpha_4^{-1}(\Omega)(ht_R(r)).
\end{eqnarray*}
This proves that $\alpha_4^{-1}(\Omega)$ is right $\hH_R$-colinear, hence
right $\hH$-colinear by \prref{prop:ff}.
By the right $\hH_R$-colinearity of $q\in
{}_{L^{op}}\mathrm{Hom}^\hH(H^{tw},A)$, the right $\hH_L$-colinearity of
$\Delta_R$, the second antipode axiom in \deref{def:hgd}(iv), 
the Hopf algebroid axiom $s_L= t_R \circ \varepsilon_R \circ s_L$ in
\deref{def:hgd}(i), the left $R$-linearity of $q$ and counitality of
$\Delta_L$, for any $h\in H$, 
\begin{eqnarray*}
(1_A\stac R h^{(1)})q(h^{(2)})
&=& q(h^{(2)})^{[0]} \stac R h^{(1)} q(h^{(2)})^{[1]}
= q({h^{(2)}}_{(2)}) \stac R h^{(1)} S({h^{(2)}}_{(1)})\\
&=& q(h_{(2)}) \stac R {h_{(1)}}^{(1)} S({h_{(1)}}^{(2)})
= q(h_{(2)}) \stac R s_L(\varepsilon_L(h_{(1)}))
= q(h) \stac R 1_H.
\end{eqnarray*}
That is, using the isomorphism \equref{eq:co_iso}, we conclude that
$\Omega^q_{A\ot_R H}(1_A \ot_R h)=q(h_{(1)})\ot_L h_{(2)}$. 
With this identity at hand, by the right $L$-linearity of $q$ and 
counitality of the coproduct in $\hH_L$, it follows that $(\alpha_4^{-1}\circ
\alpha_4)(q)=q$. On the other hand, by \equref{eq:n_11}, \equref{eq:n_13} and
the right $\hH_L$-colinearity of $\Omega_M$,
for $\Omega \in \mathrm{Nat}(U,V)$ and $m\in M$,
\begin{eqnarray*}
(\alpha_4 \circ \alpha_4^{-1})(\Omega)_M(m)
&=& {m_{[0]}}^{[0]} (A\stac L \varepsilon_L)\big(\Omega_{A\ot_R H}(1_A\stac R
  {m_{[0]}}^{[1]})\big)\stac L m_{[1]}\\
&=& (M\stac L \varepsilon_L)\big(\Omega_{M\ot_R H}({m_{[0]}}^{[0]} \stac R
{m_{[0]}}^{[1]})\big)\stac L m_{[1]}\\
&=& (M\stac L \varepsilon_L\stac L H)\big(\Omega_M(m_{[0]}) \stac L
m_{[1]}\big) \\
&=& \big( (M\stac L \varepsilon_L\stac L H)\circ (M\stac L \Delta_L)\big)
\big(\Omega_M(m)\big)
= \Omega_M(m).
\end{eqnarray*} 

It remains to check the compatibility of the constructed isomorphisms with the
connecting maps of both Morita contexts. For $p\in {}_L\mathrm{Hom}^\hH(H,A)$
and $q\in {}_{L^{op}}\mathrm{Hom}^\hH(H^{tw},A)$, $m\in M$ and $n\ot_L h\in
M^{co\hH_R}\ot_L H$, 
\begin{eqnarray*}
(\Theta^p_M\circ \Omega^q_M)(m)
&=& {m_{[0]}}^{[0]} q({m_{[0]}}^{[1]}) p(m_{[1]})
= m^{[0]} q({m^{[1]}}_{(1)}) p({m^{[1]}}_{(2)})\\
&=& m^{[0]} (q\blacklozenge' p)(m^{[1]})
= \Xi^{q\blacklozenge' p}_M(m);\\
(\Omega^q_M \circ \Theta^p_M)(n\stac L h)
&=&  n {p(h)_{[0]}}^{[0]} q \left({p(h)_{[0]}}^{[1]}\right)
\stac L p(h)_{[1]}
= np({h_{(1)}}^{(1)})q({h_{(1)}}^{(2)}) \stac L h_{(2)}\\
&=& n(p\lozenge' q)(h_{(1)}) \stac L h_{(2)}
= \Phi^{p\lozenge' q}_M(n\stac L h).
\end{eqnarray*}
This completes the proof.
\end{proof}

\prref{pr:iso_mor} justifies that \equref{eq:cat_morita} is a well defined
Morita context, for any Hopf algebroid $\hH$ and right $\hH$-comodule algebra
$A$ such that there is an algebra map from the base algebra $L$ of the
constituent left bialgebroid $\hH_L$ to the $\hH_R$-coinvariant subalgebra $B$
of $A$. 
In light of \cite[Lemmata 3.7 and 3.8]{BohmBrz:hgd.cleft}, the existence of
mutually inverse elements in the isomorphic Morita contexts in
\prref{pr:iso_mor} is equivalent to the cleft property of the $\hH$-extension
$B\subseteq A$ in the sense of \deref{def:hgd.cleft}.

\begin{example}\exlabel{ex:reg_cleft}
\prref{pr:iso_mor} implies in particular that, for any Hopf algebroid $\hH$
with constituent left bialgebroid $\hH_L$ over a base algebra $L$ and right
bialgebroid $\hH_R$ over a base algebra $R$, the $\hH$-extension
$L\subseteq H$, given by the source map $s_L$ in $\hH_L$
(equivalently, the extension $R^{op}\subseteq H$, given by the target map
$t_R$ in $\hH_R$), 
is $\hH$-cleft. Indeed, the unit of the $R$-ring $H$ is the source map $s_R$
in $\hH_R$. $H$ is an $L$-ring via the source map $s_L$ in $\hH_L$, and the 
coinvariants $L\cong \{\ s_L(l)\ |\ l\in L\ \}\equiv
\{\ t_R(r)\ |\ r\in R\ \}$ (cf. \exref{ex:hgdGal} (1))
form an $L$-subring by axiom (i) in \deref{def:hgd}. 
Mutually inverse elements in the corresponding Morita context
\equref{eq:hgd.Mor} are provided by $j$, the identity map $H$ of the total
algebra, and ${\widetilde j}$, the antipode $S$ in $\hH$, cf. axiom (iv) in
\deref{def:hgd}.  
\end{example}

A generalization of \thref{thm:pure_cleft} to arbitrary Hopf algebroids was
obtained in  \cite[Theorem 3.12]{BohmBrz:hgd.cleft}:

\begin{theorem}\thlabel{thm:nonpure.cleft}
For any Hopf algebroid $\hH$, an $\hH$-extension $B\subseteq A$ is $\hH$-cleft
if and only if the following assertions hold.

(1) The canonical map $\mathrm{can}_R$ in \equref{eq:hgd.cans} is bijective.

(2) There exists a left $B$-linear and right $\hH$-colinear isomorphism
$A\cong B \ot_L H$ (where $B \ot_L H$ is a left $B$-module via the first
factor and a right $\hH$-comodule via the second factor).
\end{theorem}

\begin{proof} (Sketch.)
Assume first that  $B\subseteq A$ is an $\hH$-cleft extension, i.e. there
exist mutually inverse elements $j\in {}_L\mathrm{Hom}^\hH(H,A)$ and 
${\widetilde j}\in {}_{L^{op}} \mathrm{Hom}^\hH(H^{tw},A)$ in the Morita
context \equref{eq:hgd.Mor}.

The inverse of $\mathrm{can}_R$ is given by 
$$
\mathrm{can}_R^{-1}(a\stac R h)=a{\widetilde j}(h_{(1)})\stac B j(h_{(2)}).
$$ 
By the $R$-$L$ bilinearity of ${\widetilde j}$ (cf. \equref{eq:bim.R-L}) and
the left $L$-linearity of $j$ (cf. \equref{eq:bim.L-R}), there is a left
$R$-module map $H\to A\ot_ L A$, $h\mapsto {\widetilde j}(h_{(1)})\ot_L
j(h_{(2)})$. Composing it with the canonical (left $R$-linear) epimorphism
$A\ot_ L A \to A\ot_ B A$, we conclude that the map $\mathrm{can}_R^{-1}$ is
well defined.

A left $B$-linear and right $\hH$-colinear isomorphism is given by the
mutually inverse maps
$$
\kappa:A\to B\stac L H,\ a\mapsto {a_{[0]}}^{[0]}{\widetilde j}({a_{[0]}}^{[1]})
\stac L a_{[1]}\qquad \textrm{and}\qquad 
\kappa^{-1}:B\stac L H \to A,\ b\stac L h \mapsto bj(h).
$$
Recall from \equref{eq:q_coinv} that $\kappa$ has the required range.

Conversely, assume that $\mathrm{can}_R$ is bijective and there exists a left
$B$-linear and right $\hH$-colinear isomorphism $\kappa:A\to B\ot_L H$. Then
mutually inverse elements $j\in {}_L\mathrm{Hom}^\hH(H,A)$ and  
${\widetilde j}\in {}_{L^{op}} \mathrm{Hom}^\hH(H^{tw},A)$ in the Morita
context \equref{eq:hgd.Mor} are given by 
$$
j:= \kappa^{-1}(1_A\stac B -)\qquad 
{\widetilde j}:= \big(A\stac B (B\stac L \varepsilon_L)\circ \kappa\big)
\big( \mathrm{can}_R^{-1}(1_A \stac R -) \big).
$$
\end{proof}

\thref{thm:hgd.cleft.str.str.thm} extends to arbitrary Hopf algebroids as
follows, see \cite[\emph{Corrigendum}]{BohmVer:Mor&cleft}. 

\begin{theorem}\thlabel{thm:nonpure.str.str.thm}
Let $\hH$ be a Hopf algebroid and $B\subseteq A$ be an $\hH$-cleft
extension. Then there is an equivalence functor $-\ot_B A:\cM_B \to \cM^\hH_A$.
\end{theorem}

\begin{proof}
Consider a cleaving map $j\in {}_L\mathrm{Hom}^\hH(H,A)$ and its inverse
${\widetilde j}\in {}_{L^{op}} \mathrm{Hom}^\hH(H^{tw},A)$ in the Morita
context \equref{eq:hgd.Mor}.

For any right $B$-module $W$, $W\ot_B A$ is an object in $\cM^\hH_A$ via the
$A$-action and the $\hH$-coactions on the second factor. The resulting functor
$-\ot_B A:\cM_B \to \cM^\hH_A$ is left adjoint of the $\hH_R$-coinvariants
functor $(-)^{co\hH_R}: \cM^\hH_A  \to \cM_B$, where the $B$-action on
$M^{co\hH_R}$ is induced by the $A$-action on $M\in \cM^\hH_A$. The counit and
the unit of the adjunction are given by
$$
\begin{array}{lll}
c_M:M^{co\hH_R}\stac B A\to M,\qquad &n \stac B a \mapsto na,\qquad
  &\textrm{for } M \in \cM^\hH_A;\\
e_W:W \to (W\stac B A)^{co\hH_R},\qquad &w \mapsto w\stac B 1_A,\qquad
  &\textrm{for } W\in \cM_B.
\end{array}
$$
The map $c_M$ is a morphism in $\cM^\hH_A$ by \coref{cor:hgd.coinv}.
We prove that $-\ot_B A:\cM_B \to \cM^\hH_A$ is an equivalence by constructing
the inverse of the above natural transformations.

The inverse of $c_M$ is given by the map
$$
c_M^{-1}(m)=m^{[0]} \mathrm{can}_R^{-1}(1_A \stac R m^{[1]})
= m^{[0]} {\widetilde j}({m^{[1]}}_{(1)}) \stac B j(m{^{[1]}}_{(2)})
= {m_{[0]}}^{[0]} {\widetilde j}({m_{[0]}}^{[1]}) \stac B j(m_{[1]}), 
$$
where we used that $\mathrm{can}_R$ is an isomorphism,
cf. \thref{thm:nonpure.cleft}. 
The third equality follows by the right $\hH_L$-colinearity of the
$\hH_R$-coaction on $M$. 
By \equref{eq:q_coinv}, the range of $c_M^{-1}$ is in $M^{co\hH_R}\ot_B A$, as
needed. For $m\in M$, and $n\ot_B a\in M^{co\hH_R}\ot_B A$,  
\begin{eqnarray*}
(c_M\circ c_M^{-1})(m)
&=& m^{[0]} {\widetilde j}({m^{[1]}}_{(1)}) j({m^{[1]}}_{(2)}) 
= m^{[0]} s_R(\varepsilon_R(m^{[1]}))
= m;\\
( c_M^{-1} \circ c_M)(n\stac B a)
&=&na^{[0]}\mathrm{can}_R^{-1}(1_A \stac R a^{[1]})
=n \, \mathrm{can}_R^{-1}(a^{[0]} \stac R a^{[1]})
= n\stac B a,
\end{eqnarray*}
where the left $A$-linearity and the explicit form of $\mathrm{can}_R$ is
used. 

The inverse of $e_W$ is given by 
$$
e_W^{-1}( \sum_i w_i\stac B a_i) 
= \sum_i w_i \, a_i^{[0]}{\widetilde j}(a_i^{[1]})j(1_H). 
$$
By \equref{eq:q_coinv}, 
$a^{[0]}{\widetilde j}(a^{[1]})$ belongs to $B$, for any $a\in A$, and by 
the $\hH_R$-colinearity of $j$ and the unitality of the right
$\hH_R$-coproduct $\Delta_R$, the element $j(1_H)$ belongs to $B$ as well.
Hence the expression of $e_W^{-1}$ is meaningful.
Obviously, $(e_W^{-1}\circ e_W)(w)=w$, for all $w\in W$. To check that
$e_W^{-1}$ is also the right inverse of $e_W$, note that for $\sum_i w_i\ot_B
a_i \in (W\ot_B A)^{co\hH_R}$ the identity 
$\sum_i w_i \ot_ B a_i^{[0]}\ot_ R a_i^{[1]}
= \sum_i w_i \ot_ B a_i \ot_ R 1_H$ holds, 
which implies 
$$
(e_W\circ e_W^{-1})(\sum_i w_i\stac B a_i)
= \sum_i w_i  a_i^{[0]}{\widetilde j}(a_i^{[1]}) j(1_H) \stac B 1_H 
= \sum_i w_i \stac B a_i^{[0]}{\widetilde j}(a_i^{[1]}) j(1_H)
= \sum_i w_i \stac B a_i.
$$  
\end{proof}

\begin{remark}\relabel{rem:fund.thm}
Since by \exref{ex:reg_cleft} the total algebra of any Hopf algebroid
$\hH$ is an $\hH$-cleft extension of the base algebra $L$,
\thref{thm:nonpure.str.str.thm} implies, in particular, that for any Hopf
algebroid $\hH$, there is an equivalence functor $-\ot_L H: \cM_L \to
\cM^\hH_H$. This yields a {\em corrected version} of the Fundamental Theorem
of Hopf modules, \cite[Theorem 4.2]{Bohm:hgdint}.

Regrettably, the proof of the journal version of \cite[Theorem
  4.2]{Bohm:hgdint} turned out to be incorrect: When checking that the
to-be-inverse of the counit of the adjunction has the required range,
some ill-defined maps are used.
Hence in the given (stronger form) \cite[Theorem 4.2]{Bohm:hgdint}
is not justified.
The journal version of \cite[Theorem 4.2]{Bohm:hgdint} becomes
justified only under some further purity assumption, see the revised arXiv
version.
\end{remark}
}
\bigskip

It remains to extend \thref{thm:cleft=cr.prod} to algebra extensions by a Hopf
algebroid. This requires first to understand the construction of a crossed
product with a bialgebroid. The following notions, introduced in
\cite[Definitions 4.1 and 4.2]{BohmBrz:hgd.cleft} extend
\deref{def:Hopf.cr.prod}. For the coproduct in a left $L$-bialgebroid $\hH$,
evaluated on an element $h$, the index notation $\Delta(h)=h_{(1)}\otimes_L
h_{(2)}$ is used.  

\begin{definition} \delabel{def:hgd.cr.prod}
A left $L$-bialgebroid $\hH$, with $L\otimes_k L^{op}$-ring structure $(H,s,t)$
and $L$-coring structure $(H,\Delta,\epsilon)$, {\em measures} an $L$-ring
$B$ with unit map $\iota:L\to B$ if there exists a $k$-module map $\cdot:
H\otimes_k B\to B$, the so called {\em measuring}, such that, for $h\in H$,
$l\in L$ and $b,b'\in B$,   

(1) $h\cdot 1_B=\iota\circ \epsilon(h) 1_B$,
 
(2) $(t(l)h)\cdot b =(h\cdot b)\iota(l)$ and $(s(l)h)\cdot b = \iota(l)
 (h\cdot b)$, 

(3) $h\cdot (bb')=(h_{(1)}\cdot b)(h_{(2)}\cdot b')$.

A $B$-valued {\em 2-cocycle} on $H$ is a $k$-module map $\sigma:
H\otimes_{L^{op}} H\to B$ (where the right and left $L^{op}$-module structures
in $H$ are given via $t$), such that, for $h,k,m\in H$ and $l\in L$,

(4) $\sigma(s(l)h,k)=\iota(l)\sigma(h,k)$ and
$\sigma(t(l)h,k)=\sigma(h,k)\iota(l)$, 

(5) $(h_{(1)}\cdot \iota(l))\sigma(h_{(2)},k)=\sigma(h,s(l)k)$,

(6) $\sigma(1_H,h)=\iota\circ \epsilon(1_H)=\sigma(h,1_H)$, 

(7) $\big(h_{(1)}\cdot\sigma(k_{(1)},m_{(1)})\big)\sigma(h_{(2)}, k_{(2)}
m_{(2)}) = \sigma(h_{(1)},k_{(1)})\sigma(h_{(2)}k_{(2)}, m)$.

The $\hH$-measured $L$-ring $B$ is a {\em $\sigma$-twisted $\hH$-module} if in
addition, for $b\in B$ and $ h,k\in H$,

(8) $1_H \cdot b =b$,

(9) $\big(h_{(1)}\cdot(k_{(1)}\cdot b)\big)\sigma(h_{(2)},k_{(2)})=
\sigma(h_{(1)},k_{(1)}) (h_{(2)}k_{(2)}\cdot b)$.

Note that condition (3) makes sense in view of (2). Conditions (7) and (9)
make sense in view of (2), (4) and (5).
\end{definition}

Similarly to \prref{prop:Cr.prod}, the reader can easily verify following
\cite[Proposition 4.3]{BohmBrz:hgd.cleft}. 

\begin{proposition}\prlabel{prop:hgd.Cr.prod}
Consider a left $L$-bialgebroid $\hH$ and an $\hH$-measured $L$-ring $B$.
Let $\sigma:H\otimes_{L^{op}} H\to B$ be a $k$-module map satisfying properties
(4) and (5) in \deref{def:hgd.cr.prod}. Take the $L$-module tensor product
$B\otimes_L H$, where $H$ is a left $L$-module via the source map.
The $k$-module $B\otimes_L H$ is a $k$-algebra with multiplication 
$$
(b\#h)(b'\#h')=b(h_{(1)}\cdot b')\sigma(h_{(2)},h'_{(1)})\#
h_{(3)}h'_{(2)}
$$
and unit element $1_B\# 1_H$, if and only if $\sigma$ is a $B$-valued
2-cocycle on $\hH$ and $B$ is a $\sigma$-twisted $\hH$-module.
This algebra is called the {\em crossed
  product} of $B$ with $\hH$ with respect to the cocycle $\sigma$. It is
denoted by $B\#_\sigma \hH$.
\end{proposition}

Similarly to the case of a crossed product with a bialgebra, 
a crossed product algebra $B\#_\sigma \hH_L$
with a left $L$-bialgebroid $\hH_L$ is  a right $\hH_L$-comodule algebra. The 
coaction is given in terms of the coproduct $\Delta_L$ in $\hH$ as $B\otimes_L
\Delta_L$. What is more, $B\subseteq B\#_\sigma \hH_L$ is an extension by
$\hH_L$.  
If $\hH_L$ is a constituent left bialgebroid in a Hopf algebroid, then 
$B\subseteq B\#_\sigma \hH_L$ is an $\hH$-extension, with respect to the
$\hH$-comodule structure on $B\#_\sigma \hH_L$ induced by the two
coproducts in $\hH$.
Analogously to the Hopf algebra case, we ask what $\hH$-extensions arise as
crossed products. Following \thref{thm:hgd.cleft=cr.prod} (extending
\thref{thm:cleft=cr.prod}) was proven in 
\cite[Theorems 4.11 and 4.12]{BohmBrz:hgd.cleft}. Before stating it, we need
to recall from \cite[Definition 4.8]{BohmBrz:hgd.cleft} what is meant by
invertibility of a $B$-valued 2-cocycle on $\hH$. In contrast to the bialgebra
case, it can not be formulated as convolution invertibility.
{\color{blue}
For its interpretation as invertibility in an appropriately chosen category,
see \cite[Definition 4.28]{Bohm:HoA}.
}

\begin{definition}
Let $\hH$ be a left $L$-bialgebroid, with $L\otimes_k L^{op}$-ring structure
$(H,s,t)$ and $L$-coring structure $(H,\Delta,\epsilon)$, and
let $\iota:L\to B$ be an $L$-ring, measured by $\hH$. A
$B$-valued 2-cocycle $\sigma$ on $\hH$ is {\em
invertible} if there exists a $k$-linear map ${\tilde\sigma}:
H\otimes_L H\to B$ (where the $L$-module structures
on $H$ are given by 
$s$), satisfying for all $h,k\in H$ and
$l\in L$, 

(1) ${\tilde \sigma}(s(l)h,k)=\iota(l) {\tilde \sigma}(h,k)$ and
${\tilde \sigma}(t(l)h,k)={\tilde\sigma}(h,k)\iota(l)$,

(2) ${\tilde \sigma}(h_{(1)},k)(h_{(2)} \cdot \iota(l))={\tilde \sigma}(h,
    t(l)k)$, 

(3) $\sigma(h_{(1)},k_{(1)})\,{\tilde \sigma}(h_{(2)},k_{(2)})=
h\cdot(k\cdot 1_B)$ and
${\tilde \sigma}(h_{(1)},k_{(1)})\,\sigma(h_{(2)},k_{(2)})=hk\cdot 1_B$.

Conditions in (3) make sense in view of (1) and (2). A map $\tilde{\sigma}$ is
called an {\em  inverse} of $\sigma$. 
\end{definition}

\begin{theorem}\thlabel{thm:hgd.cleft=cr.prod}
An algebra extension $B\subseteq A$ by a Hopf algebroid $\hH$ is a cleft
extension if and only if $A$ is isomorphic, as a left $B$-module and right
$\hH$-comodule algebra, to $B\#_\sigma \hH_L$, a crossed product of $B$ with
the constituent left bialgebroid $\hH_L$ in $\hH$ with respect to some
invertible $B$-valued 2-cocycle $\sigma$ on $\hH$.
\end{theorem}

\begin{proof}(Sketch.) 
For the two coproducts $\Delta_L$ and $\Delta_R$ in $\hH$
we use the two versions of Sweedler's index notation introduced after
\deref{def:hgd}. 
In order to show that the $\hH$-extension $B\subseteq
B\#_\sigma \hH_L$ is cleft, one constructs the inverse
of the $L$-$\hH$ bicomodule map $j:H\to B\#_\sigma \hH_L$, $h\mapsto 1_B\# h$ 
{\color{blue}
in the Morita context \equref{eq:hgd.Mor}}.
It has the explicit form
$$
{\widetilde j}(h):={\tilde \sigma}(S(h_{(1)})_{(1)},h_{(2)})\#S(h_{(1)})_{(2)}.
$$
In the converse direction, starting with an $\hH$-cleft extension $B\subseteq
A$, one constructs a measuring of $\hH_L$ on $B$ and an invertible $B$-valued
2-cocycle on $\hH_L$ in terms of the cleaving map $j$ and its inverse
${\widetilde j}$ 
{\color{blue}
in the Morita context \equref{eq:hgd.Mor}}.  
Explicitly, 
$$
h\cdot b:= j(h^{(1)})b {\widetilde j}(h^{(2)}),\qquad
\sigma(h,k):=j(h^{(1)})j(k^{(1)}){\widetilde j}(h^{(2)}k^{(2)}).
$$
The 2-cocycle $\sigma$ is proven to be invertible by constructing an inverse
$$
{\tilde \sigma}(h,k)=j(h^{(1)}k^{(1)}){\widetilde j}(k^{(2)}){\widetilde
  j}(h^{(2)}).
$$
An isomorphism $A\to B\#\hH_L$ of left $B$-modules and right $\hH$-comodule
algebras is 
given by the map $\kappa$ in the proof of \thref{thm:nonpure.cleft}.

For more details we refer to \cite[Theorems 4.11 and 4.12]{BohmBrz:hgd.cleft}.
\end{proof}

\section*{Acknowledgement}

This paper is based on two talks by the author. For possibilities to
deliver these talks she is grateful to the following people. For an
invitation to the conference ``New techniques in Hopf algebras and 
graded ring theory'' in Brussels, September 19-23, 2006, she would like to
express her gratitude to the organizers, Stef Caenepeel and Fred Van
Oystaeyen. Invitation, and warm hospitality during a visit in November 2006,
are thanked to Tomasz Brzezi\'nski and the members of Department of
Mathematics, University of Wales Swansea. Many thanks to Tomasz Brzezi\'nski
and Claudia Menini for helpful discussions. Financial support is acknowledged
to the Hungarian Scientific Research Fund OTKA T043159 and the Bolyai J\'anos
Fellowship.

\end{document}